\documentclass[11pt]{article}

\usepackage{amsmath,amssymb,amsthm,mathtools}
\usepackage{geometry}
\usepackage{enumitem}
\usepackage{hyperref}

\theoremstyle{plain}
\newtheorem{theorem}{Theorem}[section]
\newtheorem{proposition}[theorem]{Proposition}

\newtheorem{corollary}[theorem]{Corollary}

\theoremstyle{definition}
\newtheorem{definition}[theorem]{Definition}
\newtheorem{remark}[theorem]{Remark}
\newtheorem{example}[theorem]{Example}

\title{Lift Height and the Structure of Classification Failure}
\author{
Joe Gildea\\
Department of Computing Science and Mathematics,\\
School of Informatics and Creative Arts,\\
Dundalk Institute of Technology\\
\texttt{gildeajoe@gmail.com}}
\date{}

\begin{document}

\maketitle


\begin{abstract}

We introduce a framework for measuring the depth of classification
failure relative to a chosen invariant. A \emph{shadow map} records
visible structure while suppressing finer information; lift packages
recover this information in successive layers, and lift height measures
the depth through which classification ambiguity persists. Canonical
reduction removes redundant stages and yields an essential notion of
lift height. For nilpotent groups, the lower central series provides a natural
commutator lift tower, with lift height controlled by nilpotency class
and decreasing under refinement of the shadow. The associated graded
structure gives a Lie-theoretic interpretation of the successive
commutator layers.

The main examples show that classification depth is genuinely
unbounded. Over any field of characteristic different from two, and
for every positive integer, there exist nilpotent Lie algebras with the
same fixed four-dimensional abelianisation whose commutator lift
packages agree through the prescribed number of layers but differ at
the next layer. Moreover, such separation can occur at arbitrarily
large degree even when the complete underlying lower-central graded
vector space is fixed. Thus neither a fixed shadow nor the dimensions,
or even the underlying vector-space structure, of the lower-central
graded object impose a uniform bound on the depth at which new
commutator information may appear.

\end{abstract}



\section{Introduction}

Classification problems throughout algebra often begin by replacing a
complicated object with a simpler invariant: a quotient,
abelianisation, associated graded object, truncation, or other visible
shadow. When such an invariant fails to classify, the usual conclusion
is that information has been lost. The question considered in this
paper is more structural: \emph{How deep can classification failure be?}

Two objects may have the same visible invariant and remain
indistinguishable after one additional layer of structure is recovered.
They may remain indistinguishable after a second layer, or a third.
Thus classification failure need not be a single obstruction between
an invariant and a complete classification. It may persist through a
sequence of progressively finer distinctions.

We develop a framework for studying this phenomenon. A
\emph{shadow map} records visible information while suppressing finer
structure. Objects with the same shadow lie in a common shadow fibre,
where the residual classification problem is concentrated. A
\emph{lift package} records additional structure within such a fibre,
and successive lift packages form a \emph{lift tower}. The resulting
\emph{lift height} measures the depth through which hidden structure
persists relative to the chosen shadow.

This viewpoint shifts attention from the information retained by an
invariant to the ambiguity that remains after the invariant has been
fixed. The basic mechanism is

\[
\text{shadow}
\quad\longrightarrow\quad
\text{successive lift layers}
\quad\longrightarrow\quad
\text{progressive refinement of classification}.
\]

A second issue is essential. The length of a chosen lift tower need not
measure genuine classification depth, since some stages may make no new
distinctions. We therefore associate to a finite lift tower its
refinement chain of equivalence relations inside the shadow fibre and
remove repeated consecutive relations. This produces a canonical
reduced refinement chain and leads to \emph{essential lift height},
which records only genuine successive refinements.

The main phenomenon established in the paper is that essential
classification depth can be arbitrarily large. For every \(N\geq1\),
there exist nilpotent Lie algebras over a field of characteristic
different from \(2\), with the same fixed abelianisation shadow
\(F^4\), whose first \(N\) commutator lift layers agree but which are
distinguished at the next layer. Consequently, there is no uniform
finite bound on commutator classification depth even inside this single
fixed shadow fibre.

The phenomenon persists after substantially more information is fixed.
For arbitrarily large \(m\), we construct nilpotent Lie algebras
\(L_m\) and \(M_m\) of class exactly \(m\) whose associated
lower-central graded objects are isomorphic as graded vector spaces.
Their graded Lie brackets agree through degree \(m-1\), but the
degree-\(m\) commutator operation distinguishes them. Thus neither the
dimensions of the lower-central quotients nor even the complete
underlying lower-central graded vector space imposes a uniform bound
on the degree at which new commutator information may first appear.

These examples give concrete mathematical content to the notion of
lift depth. They show that classification failure can possess genuine
successive structure and that this structure cannot, in general, be
compressed into a bounded number of refinement stages determined by
the visible shadow. At the same time, lift height is fundamentally
relative: changing the shadow changes how much structure remains
hidden. For nilpotent groups of class \(c\), the commutator lift height
relative to the lower-central shadow
\(G/\gamma_{r+1}(G)\) is \(c-r\), and compatible refinement of shadows
can only decrease the remaining lift height.

The framework is related to several classical ways of organising
hidden algebraic structure. Central extensions and their classification
are naturally connected with cohomology and extension theory
\cite{Baer,Brown,Serre}; nilpotent groups are organised by the lower
central series and commutator calculus \cite{MHall,Robinson}; and their
associated graded structures connect group-theoretic and Lie-theoretic
methods, particularly in the setting of the Lazard correspondence
\cite{Lazard1954}. Related structural questions occur in coclass
theory, Baer invariants, and the classification of finite \(p\)-groups
\cite{PHALL,EickLeedhamGreen,LeedhamGreenMcKay}. More generally,
extension and obstruction phenomena are organised by homological and
derived constructions \cite{Brown,Serre,Weibel}.

Lift theory does not replace these constructions or identify them with
a single mechanism. Its point of departure is different. Rather than
fixing a particular extension, filtration, or cohomological theory, we
begin with a chosen shadow and study the classification ambiguity left
inside its fibres. The resulting lift packages and lift towers therefore
depend on the shadow itself. In the nilpotent setting the lower central
series supplies a particularly concrete realisation of this general
viewpoint, but the definitions are not tied to that realisation.

Nilpotent groups and Lie algebras provide the principal testing ground
in this paper. For suitable class-two groups, the information omitted
by abelianisation gives rise to a quadratic lift package
\((A,K,\beta,q)\), where \(A=G^{\mathrm{ab}}\),
\(K=[G,G]\), \(\beta\) is induced by the commutator, and \(q\) by
squaring. At higher nilpotency class, the successive lower-central
quotients and their commutator operations produce higher lift layers.
This identifies commutator lift depth with lower-central depth while
also providing the examples in which genuinely new information first
appears arbitrarily far along the tower.

The paper is organised as follows.
Section~\ref{sec:shadow} introduces shadow maps, shadow fibres, lift
packages, and refinement.
Section~\ref{sec:towers} develops lift towers, cumulative lift data,
and lift height.
Section~\ref{sec:universality} studies completeness, refinement,
functorial refinement, and factorisation of compatible invariants.
Section~\ref{sec:higher} develops commutator lift layers for nilpotent
groups, their associated graded interpretation, and the relation
between nilpotency class and commutator lift height.
Section~\ref{sec:relative} introduces relative lift height and shadow
refinement and shows that lift packages can contain information not
detected by lower-central quotient dimensions alone.
Section~\ref{sec:canonical-towers} develops canonical reduction and
essential lift height, constructs essential layers at arbitrarily large
depth over a fixed shadow, and proves that arbitrarily deep commutator
separation persists even when the underlying lower-central graded
vector space is fixed.
Finally, Section~\ref{sec:conclusion} summarises the framework and its
implications for the structure of classification failure.

\section{Shadow Maps, Fibres, and Lift Packages}
\label{sec:shadow}

The purpose of this section is to introduce the basic language of lift
theory. Rather than viewing the information omitted by a shadow map as
an accidental limitation of particular invariants, we interpret it as
a structural phenomenon arising from a visible projection, or
\emph{shadow}. We show that the hidden information discarded by a
shadow map is organised within shadow fibres and that recovering this
hidden structure naturally requires additional lift data, leading to
lift packages and their higher iterations. Throughout this section let
\(S:\mathcal C\to\mathcal D\) be a functor.


\subsection{Shadow maps}

The starting point of lift theory is the observation that many natural
mathematical constructions proceed by replacing a complicated object
with a simpler visible image. Such constructions preserve important
visible structure while discarding higher-order information. This
motivates the notion of a \emph{shadow map}, which formalises the
passage from an object to its visible shadow while leaving hidden
structure unresolved.

\begin{definition}
A \emph{shadow map} is a functor $S:\mathcal C\to\mathcal D$ for which there exist non-isomorphic objects
$X,Y\in\mathcal C$ satisfying $S(X)\cong S(Y)$.
\end{definition}

\begin{remark}
A shadow map retains visible or simplified structure while suppressing
higher-order information. The hidden structure omitted by the shadow
is the fundamental phenomenon studied throughout this paper.
\end{remark}

\begin{example}
Typical examples include group abelianisation
$G\mapsto G^{\mathrm{ab}}$
and associated graded constructions
$X\mapsto \operatorname{gr}(X)$.
In each case the shadow retains useful visible structure while
discarding information required for complete classification.
\end{example}

\begin{remark}
In each example, the shadow retains useful visible information while
failing to distinguish all isomorphism classes. The information omitted
by the shadow is not distributed throughout the category; it is
concentrated among objects sharing the same shadow.

This leads naturally to the notion of a shadow fibre, which isolates
the hidden structure associated with a fixed visible shadow.
\end{remark}


\subsection{Shadow Fibres}

The hidden information omitted by a shadow map is organised by its
shadow fibres. Objects sharing the same visible shadow are collected
into a common fibre, and it is within these fibres that the hidden
structure discarded by the shadow resides.

\begin{definition}
For \(A\in\mathcal D\), the \emph{shadow fibre} over \(A\) is the full
subcategory
\[
\mathcal C_A
=
\{X\in\mathcal C \mid S(X)\cong A\}.
\]
Its morphisms are all morphisms in \(\mathcal C\) between objects of
\(\mathcal C_A\).
\end{definition}

Objects in a fixed fibre share the same visible shadow data. Any
classification ambiguity therefore occurs inside the fibre.

\begin{proposition}\label{prop:classification-localisation}
Let $S:\mathcal C\to\mathcal D$ be a functor, and for \(A\in\mathcal D\) let $\mathcal C_A=\{X\in\mathcal C \mid S(X)\cong A\}$
denote the corresponding shadow fibre. Then any classification ambiguity not detected by \(S\) is contained within the fibres
\(\mathcal C_A\). Equivalently, classification relative to the shadow map \(S\) decomposes
into the classification problems inside the individual shadow fibres.
\end{proposition}

\begin{proof}
Let \(X,Y\in\mathcal C\). If $S(X)\not\cong S(Y)$, then \(X\) and \(Y\) are distinguished by the shadow map itself. Therefore, any pair of objects that remains indistinguishable after applying \(S\) must satisfy $S(X)\cong S(Y)$. Choose an object \(A\in\mathcal D\) such that $A\cong S(X)\cong S(Y)$. Then \(X\) and \(Y\) both belong to the fibre $\mathcal C_A$. Consequently, any residual classification ambiguity left unresolved by the shadow map occurs
among objects belonging to a common fibre.

Thus, once the shadow \(A\) is fixed, the remaining classification problem is entirely internal
to the fibre \(\mathcal C_A\). Hence classification relative to the shadow map \(S\) reduces to
classification within the individual fibres.
\end{proof}

Objects sharing the same shadow lie in a common shadow fibre. The
central problem of lift theory is to determine what additional
structure is required to recover the hidden information omitted by the
shadow.


\subsection{Lift packages}

The hidden information omitted by a shadow map is localised within its
individual shadow fibres. The next step is therefore to identify the
additional structure required to recover this hidden information. This
motivates the notion of a lift package.

\begin{definition}
Let $S:\mathcal C\to\mathcal D$ be a shadow map. A \emph{lift package} for \(S\) is a functor $\mathcal L:\mathcal C\to\mathcal E$
which records additional structure beyond the visible shadow and is intended to distinguish objects that are not separated by \(S\).
\end{definition}

Typically, a lift package extends the visible shadow by incorporating
operations, extension data, or higher-order interactions invisible to
\(S\). Its purpose is to recover the hidden information omitted by the
shadow and organise it into explicit mathematical structure.

\begin{definition}
A lift package is called:

\begin{enumerate}[label=(\roman*)]

\item \emph{complete} if, whenever \(X\) and \(Y\) belong to the same
shadow fibre,
\[
\mathcal L(X)\cong\mathcal L(Y)
\quad\Longrightarrow\quad
X\cong Y;
\]

\item \emph{minimal} if no strictly weaker lift package remains
complete on the same shadow fibre;

\item \emph{universal} if every compatible classification invariant
factors through it.

\end{enumerate}
\end{definition}

\begin{remark}
The purpose of a lift package is not merely to accumulate invariants.
Rather, it records the operations naturally arising from the hidden
information omitted by a shadow map.

The structure encoded by a lift package is organised into successive
lift layers, giving rise to lift towers and lift height.
\end{remark}

\subsection{Lift operations}

Having introduced lift packages abstractly, we now examine the form
that the hidden information omitted by a shadow map typically takes. A
lift package is often expressed schematically in the form
\[
\mathcal L(X)
=
(A,K,\theta_1,\theta_2,\ldots),
\]
where \(A=S(X)\) is the visible shadow, \(K\) denotes the hidden
structure omitted by the shadow, and the operations \(\theta_i\) encode
interactions invisible to the shadow map.

The first lift operations recover the most immediate hidden structure
omitted by the shadow. If these operations do not recover all of the
hidden information within a shadow fibre, further hidden structure
gives rise to higher lift operations. Iterating this process leads
naturally to lift towers and ultimately to lift height.

\subsection{Quadratic lift packages}

Let \(X\in\mathcal C\), and let \(S(X)=A\) be its shadow. The shadow
records the visible information retained by \(S\), while omitting
higher-order structure. The hidden information omitted by the shadow is
encoded in a lift package
\(\mathcal L(X)=(A,K,\theta_1,\theta_2,\ldots)\), where \(K\) denotes
hidden structure and the operations \(\theta_i\) record interactions
invisible to the shadow map. For the shadow map
\(S(G)=G^{\mathrm{ab}}\), suppose that \(G\) has nilpotency class at
most \(2\), that \(x^2\in [G,G]\) for every \(x\in G\), and that
\(c^2=1\) for every \(c\in [G,G]\). The visible shadow is
\(A=G^{\mathrm{ab}}\), while the hidden structure is
\(K=[G,G]\). In this example, the first lift operations are realised by
the commutator map \(\beta\) and a quadratic map
\(q:A\to K\). For \(a,b\in A\), choose lifts \(x,y\in G\) and define
\(\beta(a,b)=[y,x]\), \(q(a)=x^2\). Since \(G\) has nilpotency class at
most \(2\), the commutator subgroup \([G,G]\) is central, so
\(\beta(a,b)\) is independent of the chosen lifts. Moreover,
\(x^2\in [G,G]\) by assumption. If \(x'\) is another lift of \(a\),
then \(x'=xc\) for some \(c\in [G,G]\). Since \([G,G]\) is central
and \(c^2=1\), we have \((xc)^2=x^2c^2=x^2\), so \(q(a)\) is also
independent of the chosen lift. Thus the abstract lift package
\(\mathcal L(X)=(A,K,\theta_1,\theta_2,\ldots)\) takes the concrete
form \(\mathcal L(G)=(A,K,\beta,q)\), and both maps are well defined.

\begin{proposition}\label{prop:quadratic}
The operations \(\beta\) and \(q\) satisfy
\[
q(a+b)
=
q(a)+q(b)+\beta(a,b).
\]
\end{proposition}

\begin{proof}
Let \(x,y\in G\) be lifts of \(a,b\in A\). Since \(G\) has nilpotency
class at most \(2\), the commutator subgroup is central, and hence
\((xy)^2=x^2y^2[y,x]\). By definition, \(q(a+b)=(xy)^2\),
\(q(a)=x^2\), \(q(b)=y^2\), and \(\beta(a,b)=[y,x]\). Passing to
additive notation in \(K=[G,G]\), the identity becomes
\[
q(a+b)
=
q(a)+q(b)+\beta(a,b).
\]
\end{proof}

\begin{remark}
The pair \((\beta,q)\) forms the first lift layer associated with the
shadow map \(G\longmapsto G^{\mathrm{ab}}\). The key phenomenon is that
the noncommutative structure omitted by abelianisation reappears as
quadratic data on the abelian shadow.

This provides the first concrete example of a lift package and
illustrates the general principle that the hidden information omitted
by a shadow map gives rise to new algebraic structure.
\end{remark}

\subsection{General lift packages}

The quadratic lift package provides the first concrete example of the
general framework introduced above. In the notation
\(\mathcal L(X)=(A,K,\theta_1,\theta_2,\ldots)\), the operations
\(\theta_i\) are not prescribed in advance. Rather, they depend on the
chosen shadow map and the particular hidden structure under
consideration.

The guiding principle is that lift operations arise from the
information omitted by the shadow. Once a shadow map has been fixed,
objects sharing the same shadow lie in a common fibre. Any hidden
structure distinguishing such objects must therefore be encoded by
additional structure beyond the shadow itself. The operations
\(\theta_i\) are intended to capture this hidden structure.

In the quadratic setting, where the shadow is the abelianisation map
\(G\longmapsto G^{\mathrm{ab}}\), the information omitted by the
shadow is noncommutative structure. The corresponding lift operations
are the commutator map \(\beta\) and the quadratic map \(q\). Thus the
abstract lift package \((A,K,\theta_1,\theta_2,\ldots)\) takes the
concrete form \((A,K,\beta,q)\).

For other shadow maps, different hidden structures may arise and
therefore different lift operations may be required. The purpose of
the notation \(\theta_i\) is to emphasise that lift theory does not
begin with a predetermined list of operations. Instead, the relevant
operations are discovered by analysing the hidden structure associated
with a shadow fibre. The following proposition formalises the
necessity of such additional structure.

\begin{proposition}
Let \(S:\mathcal C\to\mathcal D\) be a shadow map, and let
\(\mathcal C_A\) be a shadow fibre. Suppose that \(X,Y\in\mathcal C_A\)
are non-isomorphic. Then the shadow \(A\) alone cannot distinguish
\(X\) and \(Y\). Consequently, any complete invariant on
\(\mathcal C_A\) must contain information not determined by the shadow
alone.
\end{proposition}

\begin{proof}
Since \(X,Y\in\mathcal C_A\), we have
\(S(X)\cong A\cong S(Y)\). Thus the shadow map assigns the same visible
data to \(X\) and \(Y\), and therefore does not distinguish them.
However, \(X\not\cong Y\) by assumption. Hence any complete invariant
for objects in \(\mathcal C_A\) must distinguish \(X\) from \(Y\).
Since this distinction is not detected by the shadow \(A\), it must be
made using additional data beyond the shadow.
\end{proof}

Residual lift data need not completely resolve the ambiguity within a
shadow fibre. If distinct objects remain indistinguishable after a
first lift package has been introduced, then further lift data is
required.

\begin{proposition}
Let \(S:\mathcal C\to\mathcal D\) be a shadow map, and let
\(\mathcal L_1\) be a lift package on a shadow fibre \(\mathcal C_A\).
Suppose there exist non-isomorphic objects \(X,Y\in\mathcal C_A\) such
that \(\mathcal L_1(X)\cong\mathcal L_1(Y)\). Then \(\mathcal L_1\) is
not a complete invariant on \(\mathcal C_A\). Consequently, any
complete invariant on \(\mathcal C_A\) requires additional data beyond
that recorded by \(\mathcal L_1\).
\end{proposition}

\begin{proof}
Since \(\mathcal L_1(X)\cong\mathcal L_1(Y)\), the lift package
\(\mathcal L_1\) does not distinguish \(X\) and \(Y\). However,
\(X\not\cong Y\) by assumption. Therefore \(\mathcal L_1\) cannot be a
complete invariant on \(\mathcal C_A\). Any complete invariant on
\(\mathcal C_A\) must distinguish \(X\) from \(Y\). Hence further
information is required beyond that contained in \(\mathcal L_1\).
\end{proof}

The preceding proposition motivates a comparison between different
amounts of lift data. Since the purpose of lift data is to resolve
classification ambiguity within a fixed shadow fibre, the appropriate
notion of comparison is determined by the distinctions that the data
make between objects.

\begin{definition}\label{def:lift-refinement}
Let \(S:\mathcal C\to\mathcal D\) be a shadow map, let
\(\mathcal C_A\) be a fixed shadow fibre, and let \(\mathcal L\) and
\(\mathcal M\) be lift packages on \(\mathcal C_A\). We say that
\(\mathcal M\) \emph{refines} \(\mathcal L\) if, for all
\(X,Y\in\mathcal C_A\), the implication
\(\mathcal M(X)\cong\mathcal M(Y)\Longrightarrow
\mathcal L(X)\cong\mathcal L(Y)\) holds. Thus every distinction made
by \(\mathcal L\) is also made by \(\mathcal M\).

The refinement is \emph{strict} if there exist
\(X,Y\in\mathcal C_A\) such that
\(\mathcal L(X)\cong\mathcal L(Y)\) but
\(\mathcal M(X)\not\cong\mathcal M(Y)\).
\end{definition}

\begin{remark}
Refinement therefore compares lift data by the classification
ambiguity that it leaves unresolved. A strict refinement resolves at
least one ambiguity left unresolved by the coarser lift data.
Refinement need not be strict: additional lift data may be redundant
from the point of view of classification.
\end{remark}

Successive lift layers are naturally compared through their cumulative
data. An individual higher layer should not be viewed in isolation
from the preceding layers, since the classification information
available at a given stage consists of all lift data accumulated up to
that stage.

\begin{definition}\label{def:cumulative-lift-data}
Let \(S:\mathcal C\to\mathcal D\) be a shadow map, let
\(\mathcal C_A\) be a fixed shadow fibre, and let
\(\mathcal L_1,\mathcal L_2,\ldots\) be a sequence of lift layers on
\(\mathcal C_A\). For \(n\geq 1\), the \emph{cumulative lift data
through stage \(n\)} is the combined lift package
\(\mathcal L_{\leq n}\) defined by
\(\mathcal L_{\leq n}(X)=(\mathcal L_1(X),\ldots,\mathcal L_n(X))\).
Thus two objects \(X,Y\in\mathcal C_A\) have isomorphic cumulative lift
data through stage \(n\) precisely when
\(\mathcal L_{\leq n}(X)\cong\mathcal L_{\leq n}(Y)\).
\end{definition}

The cumulative formulation makes successive refinement precise.
Passing to the next stage retains all of the preceding lift data and
adds one further layer. Consequently, later cumulative stages cannot
lose distinctions already made at earlier stages.

\begin{proposition}\label{prop:cumulative-refinement}
Let \(\mathcal L_1,\mathcal L_2,\ldots\) be a sequence of lift layers
on a fixed shadow fibre \(\mathcal C_A\). For every \(n\geq 1\), the
cumulative lift package \(\mathcal L_{\leq n+1}\) refines
\(\mathcal L_{\leq n}\). Equivalently, for all
\(X,Y\in\mathcal C_A\), the implication
\(\mathcal L_{\leq n+1}(X)\cong\mathcal L_{\leq n+1}(Y)
\Longrightarrow
\mathcal L_{\leq n}(X)\cong\mathcal L_{\leq n}(Y)\) holds.
\end{proposition}

\begin{proof}
Suppose that
\(\mathcal L_{\leq n+1}(X)\cong\mathcal L_{\leq n+1}(Y)\).
By definition, the cumulative data through stage \(n+1\) consists of
the first \(n+1\) lift layers. Restricting this isomorphism to the
first \(n\) components gives
\(\mathcal L_{\leq n}(X)\cong\mathcal L_{\leq n}(Y)\).
Hence \(\mathcal L_{\leq n+1}\) refines
\(\mathcal L_{\leq n}\).
\end{proof}

\begin{remark}
The refinement in Proposition~\ref{prop:cumulative-refinement} need
not be strict. The passage from stage \(n\) to stage \(n+1\) is strict
precisely when there exist \(X,Y\in\mathcal C_A\) such that
\(\mathcal L_{\leq n}(X)\cong\mathcal L_{\leq n}(Y)\) but
\(\mathcal L_{\leq n+1}(X)\not\cong
\mathcal L_{\leq n+1}(Y)\). In this case the \((n+1)\)-st layer
resolves a classification ambiguity left unresolved by the preceding
layers. If no such pair exists, the additional stage is redundant from
the point of view of classification.
\end{remark}

These results provide a precise justification for the passage from
lift packages to lift layers. A shadow map may omit hidden structure;
a first lift package may recover part, but not all, of this hidden
information; and further lift data may then be required. The
cumulative packages \(\mathcal L_{\leq 1},
\mathcal L_{\leq 2},\ldots\) form a successive refinement chain, while
strict refinement identifies those stages at which genuinely new
classification distinctions are introduced. This iterative structure
leads naturally to lift towers and ultimately to the notion of lift
height developed later in the paper.


\section{Lift Towers and Lift Height}
\label{sec:towers}

The previous section introduced lift packages as carriers of hidden
structure omitted by a shadow map. In practice, however, a single lift
package need not recover all of this hidden structure at once.
Additional hidden interactions may remain after the first lift layer is
taken into account.

This suggests a hierarchical viewpoint. Rather than viewing lift
structure as a single object, we organise it into successive layers,
each recovering hidden structure left unresolved by the preceding
ones. This leads naturally to the notions of lift towers and lift
height.

\subsection{Lift layers}

The first step in constructing a lift tower is to isolate the
individual layers of hidden structure associated with a shadow map. Let
$S:\mathcal C\to\mathcal D$ be a shadow map, and fix a fibre
$\mathcal C_A$. The first lift operations determine an initial lift
layer $\mathcal L_1$. In general, the data contained in
\(\mathcal L_1\) does not recover all of the hidden structure within
the fibre. Additional hidden structure then gives rise to further lift
operations defining higher layers $\mathcal L_2,\mathcal L_3,\dots$.

\begin{remark}
Each successive lift layer recovers hidden structure left unresolved by
the preceding layers. Higher lift layers arise because lower-order lift
data does not capture all of the hidden interactions associated with a
shadow fibre.
\end{remark}

The collection of all lift layers records the successive stages through
which hidden structure is recovered beyond the original shadow. This
hierarchy is formalised by the notion of a lift tower.

\subsection{Lift towers}

Lift layers become most useful when viewed collectively rather than in
isolation. The successive layers arising from the hidden structure
omitted by a shadow map naturally assemble into a hierarchical
structure recording increasing levels of hidden information.

\begin{definition}
A \emph{lift tower} is a hierarchy of lift layers
$\mathcal L_1,\mathcal L_2,\dots$
organising successive stages of hidden structure.
\end{definition}

More precisely, \(\mathcal L_1\) refines the original shadow map, and
for each \(n\ge 1\), the layer \(\mathcal L_{n+1}\) records information
not detected by $(\mathcal L_1,\dots,\mathcal L_n)$. Consequently, the hierarchy encodes progressively finer layers of
hidden structure.

\begin{proposition}
\label{prop:lift-package-tower}
Every complete lift package determines a lift tower after choosing an
ordering of its lift data by successive refinement. Conversely, every lift tower determines a lift package by combining the
data of all its lift layers. If the tower resolves all remaining
ambiguity within the shadow fibre, then the resulting lift package is
complete.
\end{proposition}

\begin{proof}
Let \(\mathcal L:\mathcal C_A\to\mathcal E\) be a complete lift package on a fixed shadow fibre \(\mathcal C_A\). By completeness, the isomorphism class of \(\mathcal L(X)\) determines the isomorphism class
of \(X\) within \(\mathcal C_A\). Choose a decomposition of the lift data carried by \(\mathcal L\) into successive layers \(\mathcal L_1,\mathcal L_2,\ldots\), where \(\mathcal L_1\) records the first refinement beyond the shadow and
\(\mathcal L_{n+1}\) records the additional information not already captured by \((\mathcal L_1,\ldots,\mathcal L_n)\). Then the sequence \(\mathcal L_1,\mathcal L_2,\ldots\)
forms a lift tower. Since the combined lift layers recover the original package \(\mathcal L\), the resulting tower contains the same information as \(\mathcal L\).

Conversely, suppose a lift tower \(\mathcal L_1,\mathcal L_2,\ldots\) is given. Define its associated total lift package by $\mathcal L_{\mathrm{tot}}(X)=
(\mathcal L_1(X),\mathcal L_2(X),\ldots)$. This package records simultaneously the information contained in every lift layer and therefore determines the original tower. If the tower
resolves all remaining ambiguity within the shadow fibre, then $\mathcal L_{\mathrm{tot}}(X)\cong \mathcal L_{\mathrm{tot}}(Y)$ implies \(X\cong Y\). Hence \(\mathcal L_{\mathrm{tot}}\) 
is a complete lift package. Therefore a complete lift package may be organised into successive lift
layers, while a lift tower may be assembled into a lift package by
collecting the data of all its layers. Whenever the tower resolves all
remaining ambiguity, the resulting lift package is complete.
\end{proof}

This observation explains why lift towers provide a natural framework
for organising hidden structure. Rather than treating lift structure
as a single invariant, one separates it into successive layers
according to the hidden structure recovered at each stage.

\begin{remark}
Lift towers behave formally like obstruction towers: each layer records
the hidden structure not captured by the preceding layers. The
hierarchy reflects the fact that the hidden structure associated with a
shadow fibre organises into successive obstruction layers rather than
appearing randomly.
\end{remark}

Lift towers therefore organise hidden structure according to depth.
This makes it possible to measure lift height by asking how many layers
are required to recover all of the hidden structure associated with the
original shadow.


\subsection{Lift height}

Lift towers describe how hidden structure is distributed across
successive lift layers. This naturally leads to a numerical measure of
how deeply hidden structure is organised: the number of lift layers
required to recover all of the hidden structure associated with a
shadow fibre.

The resulting invariant is called the \emph{lift height}.

\begin{definition}
The \emph{lift height} of an object is the minimal number of lift
layers required to distinguish it within its shadow fibre.
\end{definition}

Equivalently, if \(\mathcal L_1,\mathcal L_2,\dots\) is the associated
lift tower, then
\[
h(X)
=
\min\left\{
n\geq 1
\;\middle|\;
(\mathcal L_1,\dots,\mathcal L_n)
\text{ determines }X
\right\}.
\]

\begin{proposition}
\label{prop:first-layer-height-one}
If the first lift layer already classifies objects in a shadow fibre,
then the associated lift height is \(1\).
\end{proposition}

\begin{proof}
Let $\mathcal C_A$ be a shadow fibre, and let $\mathcal L_1,\mathcal L_2,\ldots$ be the associated lift tower.
Assume that the first lift layer already classifies objects in \(\mathcal C_A\). Thus, for any objects $X,Y\in\mathcal C_A$,
the implication
\[
\mathcal L_1(X)\cong\mathcal L_1(Y)
\quad\Longrightarrow\quad
X\cong Y
\]
holds. By definition, the lift height \(h(X)\) is the minimal number of lift
layers required to determine the isomorphism class of \(X\) within its
shadow fibre. Since the first lift layer already determines the
isomorphism class of every object in \(\mathcal C_A\), one layer
suffices. Therefore the minimal such number is \(1\), and hence  $h(X)=1$.
\end{proof}

\begin{remark}
Lift height measures the depth of structure lying beyond the visible
shadow. Height-one phenomena occur when the first lift layer already
recovers all of the hidden structure associated with a shadow fibre,
whereas higher lift height indicates the presence of genuinely
higher-order obstruction structure.
\end{remark}

\[
\text{shadow}
\;\longrightarrow\;
\text{lift layers}
\;\longrightarrow\;
\text{hidden structure}.
\]

\begin{remark}
Lift towers organise hidden structure into successive layers of
refinement. The first lift layer records the most immediate information
omitted by the shadow map, while higher lift layers recover hidden
structure not detected by lower layers.

Viewed in this way, lift height records how many stages are required to
recover all of the hidden structure associated with a shadow.
\end{remark}

\subsection{Lift towers and fibre refinement}

Shadow fibres and lift towers organise information in different ways.
A shadow map decomposes objects into fibres \(\mathcal C_A\), whose
objects share the same visible shadow. These fibres identify where
hidden structure is concentrated.

Once a shadow fibre has been fixed, a lift tower does not replace the
fibre decomposition. Rather, it organises the hidden structure within
the fibre into successive lift layers. Each finite stage of the tower
determines a relation on objects in the fibre: two objects are
equivalent at stage \(n\) if they cannot be distinguished by the first
\(n\) lift layers. More precisely, for a lift tower
\[
\mathcal L_1,\mathcal L_2,\ldots
\]
on a fibre \(\mathcal C_A\), define a relation \(\sim_n\) on objects of
\(\mathcal C_A\) by
\[
X\sim_n Y
\quad\Longleftrightarrow\quad
(\mathcal L_1(X),\ldots,\mathcal L_n(X))
\cong
(\mathcal L_1(Y),\ldots,\mathcal L_n(Y)).
\]
Thus \(X\sim_n Y\) means that \(X\) and \(Y\) remain
indistinguishable after the first \(n\) lift layers.

\begin{proposition}
Let \(\mathcal L_1,\mathcal L_2,\ldots\) be a lift tower on a shadow
fibre \(\mathcal C_A\). Then the relations \(\sim_n\) form a
successively refining hierarchy:
\[
X\sim_{n+1}Y
\quad\Longrightarrow\quad
X\sim_n Y.
\]
Consequently, each lift layer can only refine, and never coarsen, the
equivalence relation induced by the preceding layers.
\end{proposition}

\begin{proof}
Suppose that \(X\sim_{n+1}Y\). By definition,
\[
(\mathcal L_1(X),\ldots,\mathcal L_n(X),\mathcal L_{n+1}(X))
\cong
(\mathcal L_1(Y),\ldots,\mathcal L_n(Y),\mathcal L_{n+1}(Y)).
\]
Forgetting the final component \(\mathcal L_{n+1}\), we obtain
\[
(\mathcal L_1(X),\ldots,\mathcal L_n(X))
\cong
(\mathcal L_1(Y),\ldots,\mathcal L_n(Y)).
\]
Hence \(X\sim_n Y\). Therefore the equivalence classes determined by
stage \(n+1\) are contained in those determined by stage \(n\), so each
new lift layer refines the equivalence relation induced by the previous
layers.
\end{proof}

This result makes precise the relationship between fibres and towers.
The shadow map first gives a coarse decomposition into fibres. Inside a
fixed fibre, the lift tower produces a sequence of finer and finer
distinctions:
\[
\mathcal C_A/\!\sim_1
\;\leftarrow\;
\mathcal C_A/\!\sim_2
\;\leftarrow\;
\mathcal C_A/\!\sim_3
\;\leftarrow\;
\cdots .
\]
The arrows indicate that each later stage refines the equivalence
classes of the earlier stages. Thus the shadow performs a coarse
classification by grouping objects according to visible data, while the
lift tower performs progressively finer classification within a fixed
shadow fibre. In this sense, shadow fibres identify where hidden
structure resides, whereas lift towers describe how that hidden
structure is progressively recovered.

The refinement process may be viewed schematically as
\[
\mathcal C
\;\xrightarrow{\;S\;}\;
\{\text{shadow fibres}\}
\;\xrightarrow{\;\mathcal L_1\;}\;
\{\text{refined fibres}\}
\;\xrightarrow{\;\mathcal L_2\;}\;
\{\text{further refined fibres}\}
\;\xrightarrow{\;\mathcal L_3\;}\;
\cdots .
\]

The shadow map first decomposes the category according to visible
information. Successive lift layers then refine the resulting fibres by
distinguishing objects that remain indistinguishable at earlier stages.
Thus the hidden information omitted by the shadow is progressively
recovered through the hierarchy of lift layers.

\section{Universality and Reconstruction of Lift Packages}
\label{sec:universality}

The preceding sections introduced shadow maps, lift packages, lift
towers, and lift height as tools for organising the hidden structure
omitted by a shadow map. The purpose of the present section is to
examine the role played by lift packages in classification,
refinement, and the organisation of invariants associated with a fixed
shadow fibre.

Three related principles are developed. First, a complete lift package
determines the isomorphism class of an object within its shadow fibre.
Second, completeness is preserved under refinement and therefore
depends only on the classification information carried by a lift
package. Third, at the functorial level, lift packages organise
invariants that factor through their data, and this organisation is
preserved under further functorial constructions and functorial
refinement.

Together, these principles describe the position of lift packages
between shadows and the objects that realise them. A shadow records
the visible structure of an object, while a lift package records
additional structure needed to distinguish objects within the
corresponding shadow fibre and, when appropriate, to organise
functorial invariants of those objects.


\subsection{Complete lift packages and reconstruction}

The purpose of a complete lift package is not necessarily to provide an
explicit construction of an object. Rather, its primary role is to
recover enough of the hidden structure associated with a fixed shadow
fibre to determine objects up to isomorphism.

Thus reconstruction in lift theory may be understood in a
classification sense. If two objects in the same shadow fibre have
isomorphic complete lift packages, then they are themselves
isomorphic. Consequently, the complete lift-package data determines
the isomorphism class of an object within that fibre.

The following proposition records this general reconstruction
principle.

\begin{proposition}
\label{prop:complete-reconstruction}
Let \(S:\mathcal C\to\mathcal D\) be a shadow map, and let
\(\mathcal C_A\) be a shadow fibre. If \(\mathcal L\) is a complete
lift package on \(\mathcal C_A\), then \(\mathcal L(X)\) determines
the isomorphism class of \(X\) inside \(\mathcal C_A\).
\end{proposition}

\begin{proof}
Let \(X,Y\in\mathcal C_A\). Since \(\mathcal L\) is complete,
\(\mathcal L(X)\cong\mathcal L(Y)\) implies \(X\cong Y\).
Thus no two non-isomorphic objects in the fibre have isomorphic
complete lift packages. Hence the data \(\mathcal L(X)\) determines
the isomorphism class of \(X\) within \(\mathcal C_A\).
\end{proof}

This principle is intrinsic to the definition of completeness. It does
not assert that an object can necessarily be constructed explicitly
from its lift-package data. Rather, it says that a complete lift
package contains sufficient information to distinguish isomorphism
classes within a fixed shadow fibre.

This distinction is useful. Explicit reconstruction is an additional
property that may depend on the category and the particular form of
the lift package, whereas classification within a shadow fibre is
precisely the role encoded by completeness.


\subsection{Refinement and completeness}

The refinement relation compares lift packages according to the
classification information they retain within a fixed shadow fibre.
A lift package that refines another makes at least all of the
distinctions made by the coarser package. It is therefore natural to
ask how completeness behaves under refinement.

We first isolate the case in which two lift packages carry exactly the
same classification information.

\begin{definition}
\label{def:refinement-equivalent}
Let \(\mathcal L\) and \(\mathcal M\) be lift packages on a fixed
shadow fibre \(\mathcal C_A\). We say that \(\mathcal L\) and
\(\mathcal M\) are \emph{refinement-equivalent} if each refines the
other. Equivalently, for all \(X,Y\in\mathcal C_A\),
\(\mathcal L(X)\cong\mathcal L(Y)\) if and only if
\(\mathcal M(X)\cong\mathcal M(Y)\).
Thus refinement-equivalent lift packages make exactly the same
distinctions between isomorphism classes within the shadow fibre.
\end{definition}

Completeness is monotone with respect to refinement.

\begin{proposition}
\label{prop:refinement-preserves-completeness}
Let \(\mathcal L\) and \(\mathcal M\) be lift packages on a fixed
shadow fibre \(\mathcal C_A\). Suppose that \(\mathcal M\) refines
\(\mathcal L\). If \(\mathcal L\) is complete, then
\(\mathcal M\) is complete.
\end{proposition}

\begin{proof}
Suppose that \(\mathcal M(X)\cong\mathcal M(Y)\). Since
\(\mathcal M\) refines \(\mathcal L\), it follows that
\(\mathcal L(X)\cong\mathcal L(Y)\). Since \(\mathcal L\) is complete,
this implies \(X\cong Y\). Therefore \(\mathcal M\) is complete.
\end{proof}

As an immediate consequence, completeness is invariant under
refinement-equivalence.

\begin{corollary}
\label{cor:refinement-equivalent-completeness}
Let \(\mathcal L\) and \(\mathcal M\) be refinement-equivalent lift
packages on a fixed shadow fibre \(\mathcal C_A\). Then
\(\mathcal L\) is complete if and only if \(\mathcal M\) is complete.
\end{corollary}

\begin{proof}
Suppose first that \(\mathcal L\) is complete. Since \(\mathcal M\)
refines \(\mathcal L\),
Proposition~\ref{prop:refinement-preserves-completeness} implies that
\(\mathcal M\) is complete. Conversely, if \(\mathcal M\) is complete,
then \(\mathcal L\) refines \(\mathcal M\), so the same proposition
implies that \(\mathcal L\) is complete.
\end{proof}

Completeness therefore depends only on the classification information
carried by a lift package, rather than on a particular presentation of
that information. In particular, replacing a lift package by a
refinement-equivalent package does not alter whether the hidden
structure of the shadow fibre has been recovered up to isomorphism.


\subsection{Universality and factorisation}

The preceding discussion concerns the classification information
carried by lift packages. There is also a stronger, functorial level
at which one lift package may determine another. This distinction is
important when lift packages are used not only to distinguish objects
but also to organise invariants and their behaviour on morphisms.

\begin{definition}
\label{def:functorial-refinement}
Let \(\mathcal L_1:\mathcal C_A\to\mathcal E_1\) and
\(\mathcal L_2:\mathcal C_A\to\mathcal E_2\) be lift packages.
We say that \(\mathcal L_2\) \emph{functorially refines}
\(\mathcal L_1\) if there exists a functor
\(R:\mathcal E_2\to\mathcal E_1\) together with a natural isomorphism
\(\mathcal L_1\cong R\circ\mathcal L_2\).
\end{definition}

Functorial refinement is stronger than the refinement relation defined
earlier. Indeed, if the finer lift package cannot distinguish two
objects, then neither can any lift package obtained functorially from
it.

\begin{proposition}
\label{prop:functorial-refinement-implies-refinement}
Let \(\mathcal L_1:\mathcal C_A\to\mathcal E_1\) and
\(\mathcal L_2:\mathcal C_A\to\mathcal E_2\) be lift packages.
If \(\mathcal L_2\) functorially refines \(\mathcal L_1\), then
\(\mathcal L_2\) refines \(\mathcal L_1\).
\end{proposition}

\begin{proof}
Suppose that \(\mathcal L_2(X)\cong\mathcal L_2(Y)\). Since
\(\mathcal L_2\) functorially refines \(\mathcal L_1\), there exists
a functor \(R:\mathcal E_2\to\mathcal E_1\) together with a natural
isomorphism \(\mathcal L_1\cong R\circ\mathcal L_2\).
Applying \(R\) to an isomorphism
\(\mathcal L_2(X)\cong\mathcal L_2(Y)\) gives
\(R(\mathcal L_2(X))\cong R(\mathcal L_2(Y))\).
The natural isomorphism
\(\mathcal L_1\cong R\circ\mathcal L_2\) therefore gives
\(\mathcal L_1(X)\cong\mathcal L_1(Y)\). Hence
\(\mathcal L_2\) refines \(\mathcal L_1\).
\end{proof}

Thus functorial refinement implies refinement at the classification
level. In particular, by
Proposition~\ref{prop:refinement-preserves-completeness}, if
\(\mathcal L_1\) is complete and \(\mathcal L_2\) functorially
refines \(\mathcal L_1\), then \(\mathcal L_2\) is also complete.

Functorial refinement concerns the relation between two lift packages.
We now use the same factorisation principle to describe invariants
whose information is controlled by a given lift package.

\begin{definition}
\label{def:compatible-invariant}
Let \(\mathcal L:\mathcal C_A\to\mathcal E\) be a lift package and let
\(F:\mathcal C_A\to\mathcal F\) be a functorial invariant. We say that
\(F\) is \emph{compatible with} \(\mathcal L\) if there exists a
functor \(\overline F:\mathcal E\to\mathcal F\) together with a natural
isomorphism \(F\cong\overline F\circ\mathcal L\).
\end{definition}

Thus compatibility means that the invariant \(F\) can be recovered
functorially from the lift package. In particular, its dependence on
both objects and morphisms is mediated by the lift data rather than
merely by the isomorphism classes of the objects involved. This
morphism-level condition is stronger than merely requiring that
isomorphic lift-package values give isomorphic values of \(F\).

Compatibility is preserved under further functorial constructions.

\begin{proposition}
\label{prop:compatibility-postcomposition}
Let \(\mathcal L:\mathcal C_A\to\mathcal E\) be a lift package and let
\(F:\mathcal C_A\to\mathcal F\) be compatible with \(\mathcal L\).
If \(H:\mathcal F\to\mathcal G\) is a functor, then
\(H\circ F\) is compatible with \(\mathcal L\).
\end{proposition}

\begin{proof}
Since \(F\) is compatible with \(\mathcal L\), there exists a functor
\(\overline F:\mathcal E\to\mathcal F\) together with a natural
isomorphism \(F\cong\overline F\circ\mathcal L\).
Postcomposing this natural isomorphism with \(H\) gives a natural
isomorphism
\(H\circ F\cong H\circ\overline F\circ\mathcal L\).
Since \(H\circ\overline F:\mathcal E\to\mathcal G\) is a functor,
\(H\circ F\) factors through \(\mathcal L\). Hence
\(H\circ F\) is compatible with \(\mathcal L\).
\end{proof}

The invariants controlled by a lift package also behave monotonically
under functorial refinement. Passing to finer functorial lift data
cannot destroy an invariant that was already determined by coarser
lift data.

\begin{proposition}
\label{prop:factorisation-refinement}
Let \(\mathcal L_1:\mathcal C_A\to\mathcal E_1\) and
\(\mathcal L_2:\mathcal C_A\to\mathcal E_2\) be lift packages.
Suppose that \(\mathcal L_2\) functorially refines
\(\mathcal L_1\). If a functorial invariant
\(F:\mathcal C_A\to\mathcal F\) is compatible with
\(\mathcal L_1\), then \(F\) is compatible with
\(\mathcal L_2\).
\end{proposition}

\begin{proof}
Since \(\mathcal L_2\) functorially refines \(\mathcal L_1\), there
exists a functor \(R:\mathcal E_2\to\mathcal E_1\) together with a
natural isomorphism \(\mathcal L_1\cong R\circ\mathcal L_2\).
Since \(F\) is compatible with \(\mathcal L_1\), there exists a functor
\(\overline F:\mathcal E_1\to\mathcal F\) together with a natural
isomorphism \(F\cong\overline F\circ\mathcal L_1\).
Postcomposing the first natural isomorphism with \(\overline F\) and
combining it with the second gives
\(F\cong\overline F\circ\mathcal L_1\cong
\overline F\circ R\circ\mathcal L_2\).
Since \(\overline F\circ R:\mathcal E_2\to\mathcal F\) is a functor,
\(F\) factors through \(\mathcal L_2\). Hence \(F\) is compatible with
\(\mathcal L_2\).
\end{proof}

These results give a precise sense in which lift packages act as
organising objects for hidden structure. At the classification level,
refinement compares the distinctions made between objects in a shadow
fibre, and completeness is preserved as that information is refined.
At the functorial level, compatibility identifies invariants whose
information factors through a lift package, while functorial
refinement preserves all invariants already organised by coarser lift
data.

The two levels are related but distinct. Ordinary refinement concerns
classification information and is formulated in terms of isomorphism
classes. Functorial refinement additionally records how one package is
obtained from another through a functor and therefore controls
morphism-level information. The implication from functorial refinement
to ordinary refinement connects these two levels without identifying
them.

In this sense, universality in lift theory is expressed through
factorisation. A lift package provides a common structure through
which compatible invariants may be expressed, and a functorially finer
package retains every such invariant available at a coarser stage.
This complements the reconstruction role of complete lift packages:
completeness measures whether enough hidden information has been
recovered to classify objects, while factorisation describes how that
information organises further invariants.

\section{Higher Lift Layers and Nilpotency Structure}
\label{sec:higher}

The purpose of this section is to explain how higher nilpotency
naturally generates higher lift layers. The successive commutator
layers arising from the lower central series admit a classical
interpretation through the associated graded Lie structure and give
rise to a natural notion of commutator lift depth.

The resulting picture reveals a close relationship between nilpotency
and the complexity of the commutator lift tower associated with the
abelianisation shadow: deeper lower-central structure produces
correspondingly deeper layers of hidden structure. We shall also see
that lift height is fundamentally a shadow-dependent notion, so this
correspondence should be viewed as a feature of the commutator-lift
setting rather than a general property of all lift towers.


\subsection{The associated graded viewpoint}

The higher commutator layers arising from the lower central series are
closely related to the classical Lie-theoretic structure associated
with a nilpotent group. Rather than appearing as unrelated pieces of
additional data, the successive layers fit naturally into the
associated graded Lie ring of the lower central series.

\begin{theorem}
\label{thm:lazard}
Let \(G\) be a nilpotent group of class \(c\), with lower central series
\(\gamma_1(G)=G\) and
\(\gamma_{i+1}(G)=[G,\gamma_i(G)]\). Define
\(A=G^{\mathrm{ab}}\) and
\(K_i=\gamma_i(G)/\gamma_{i+1}(G)\) for \(2\leq i\leq c\).
Then the layers \(A,K_2,\ldots,K_c\) are precisely the nontrivial
homogeneous components of the associated graded Lie ring
\[
\operatorname{gr}(G)
=
\bigoplus_{i\geq1}\gamma_i(G)/\gamma_{i+1}(G).
\]
Moreover, the commutator operations between the lift layers agree with
the corresponding graded Lie brackets.
\end{theorem}

\begin{proof}
Since \(G\) has nilpotency class \(c\), we have
\(\gamma_{c+1}(G)=1\), so the homogeneous components of
\(\operatorname{gr}(G)\) vanish in degrees greater than \(c\).
For all \(i,j\geq1\), the lower central series satisfies
\[
[\gamma_i(G),\gamma_j(G)]
\subseteq
\gamma_{i+j}(G).
\]
Hence the group commutator induces maps
\[
[-,-]:
\gamma_i(G)/\gamma_{i+1}(G)
\times
\gamma_j(G)/\gamma_{j+1}(G)
\longrightarrow
\gamma_{i+j}(G)/\gamma_{i+j+1}(G),
\]
given by
\[
[x\gamma_{i+1}(G),y\gamma_{j+1}(G)]
=
[x,y]\gamma_{i+j+1}(G).
\]
These maps are well defined: changing the first representative
introduces commutators lying in
\[
[\gamma_{i+1}(G),\gamma_j(G)]
\subseteq
\gamma_{i+j+1}(G),
\]
while changing the second introduces commutators lying in
\[
[\gamma_i(G),\gamma_{j+1}(G)]
\subseteq
\gamma_{i+j+1}(G).
\]

The standard commutator identities imply that the induced bracket is
bilinear and alternating and satisfies the Jacobi identity. Thus
\[
\operatorname{gr}(G)
=
\bigoplus_{i\geq1}\gamma_i(G)/\gamma_{i+1}(G)
\]
is a graded Lie ring.

Its degree-one component is
\[
\gamma_1(G)/\gamma_2(G)
=
G/[G,G]
=
G^{\mathrm{ab}}
=
A,
\]
while for \(2\leq i\leq c\) its degree-\(i\) component is
\[
\gamma_i(G)/\gamma_{i+1}(G)
=
K_i.
\]
Since the components in degrees greater than \(c\) vanish, the layers
\(A,K_2,\ldots,K_c\) are precisely the nontrivial homogeneous
components of \(\operatorname{gr}(G)\). The operations between these
layers induced by group commutators are exactly the corresponding
graded Lie brackets. Therefore the higher commutator layers
\(A,K_2,\ldots,K_c\), together with their induced commutator
operations, form the associated graded Lie structure
\(\operatorname{gr}(G)\).
\end{proof}

\begin{remark}
The theorem places the commutator lift layers within a familiar
mathematical setting. They are not additional structures imposed on a
nilpotent group. Rather, they coincide with the homogeneous components
and graded commutator operations arising from the lower central series.
Thus the higher commutator lift package may be viewed as a
lift-theoretic organisation of the associated graded Lie structure.

The significance of the theorem is therefore not the construction of a
new Lie-theoretic invariant. Instead, it shows that the lift-theoretic
framework recovers a classical structure in a setting where the hidden
commutator structure is already well understood. This provides evidence
that lift packages capture genuine algebraic structure rather than
formal bookkeeping.
\end{remark}


\subsection{Nilpotency and lift depth}

The associated graded description suggests a close relationship
between the lower central series of a nilpotent group and the
successive layers of its commutator lift tower. Each nontrivial
lower-central quotient contributes a further layer of commutator
structure beyond the abelianisation shadow.

More generally, each nontrivial stage of the lower central series
contributes an additional layer of hidden commutator structure,
requiring a corresponding lift layer in the commutator lift tower.
This relationship is made precise by the following theorem.


\begin{theorem}
\label{thm:nilpotency-lift-height}
Let \(G\) be a nilpotent group of class \(c\). Consider the commutator
lift tower associated to the abelianisation shadow
\(G\longmapsto G^{\mathrm{ab}}\), whose \(i\)-th layer is
\[
\mathcal L_i^{\mathrm{com}}(G)
=
\gamma_{i+1}(G)/\gamma_{i+2}(G).
\]
Then \(h_{\mathrm{com}}(G)=c-1\).
\end{theorem}

\begin{proof}
Since \(G\) has nilpotency class \(c\),
\(\gamma_c(G)\neq 1\), \(\gamma_{c+1}(G)=1\). The commutator lift tower
beyond the abelianisation shadow consists of the successive layers
\[
\gamma_2(G)/\gamma_3(G),
\gamma_3(G)/\gamma_4(G),
\dots,
\gamma_c(G)/\gamma_{c+1}(G).
\]

There are exactly \(c-1\) such layers. Each of these layers is
nontrivial. Indeed, if
\(\gamma_i(G)/\gamma_{i+1}(G)\) were trivial for some
\(2\leq i\leq c\), then
\(\gamma_i(G)=\gamma_{i+1}(G)\), so the lower central series would
stabilise from that point onward. Since \(\gamma_{c+1}(G)=1\), this
would force \(\gamma_i(G)=1\), contradicting
\(\gamma_c(G)\neq 1\). The final layer is
\(\gamma_c(G)/\gamma_{c+1}(G)=\gamma_c(G)\), which is nonzero, while
the next layer vanishes because \(\gamma_{c+1}(G)=1\). Thus the
commutator lift tower has precisely \(c-1\) nontrivial stages. By
definition of commutator lift height,
\(h_{\mathrm{com}}(G)=c-1\).
\end{proof}

\begin{remark}
The theorem identifies commutator lift height with ordinary nilpotency
class, up to the shift introduced by abelianisation. In particular,
each nontrivial lower-central quotient contributes exactly one
additional layer of hidden commutator structure, and hence exactly one
additional commutator lift layer.
\end{remark}

\begin{corollary}
\label{cor:unitriangular-height}
Let \(F\) be a field and let \(n\geq 2\). For the unitriangular group
\(UT_n(F)\),
\[
h_{\mathrm{com}}(UT_n(F))=n-2.
\]
\end{corollary}

\begin{proof}
The group \(UT_n(F)\) has nilpotency class \(n-1\). Applying
Theorem~\ref{thm:nilpotency-lift-height} gives
\[
h_{\mathrm{com}}(UT_n(F))
=
(n-1)-1
=
n-2.
\]
\end{proof}

The preceding theorem identifies commutator lift height with the depth
of the lower central series relative to the abelianisation shadow.
This gives a concrete interpretation of lift height in terms of
commutator structure. The dependence on the abelianisation shadow is
important, however: lift height is not an intrinsic invariant of the
underlying object.

\begin{proposition}
Lift height is a relative notion depending on the chosen shadow.
In particular, two shadow maps on the same category can assign
different lift heights to the same object.
\end{proposition}

\begin{proof}
Consider the identity shadow \(S_1=\operatorname{id}_{\mathcal C}\).
Since this shadow retains the object itself, each shadow fibre contains
only a single isomorphism class. Hence there is no residual
classification ambiguity to resolve, and every object has lift height
\(0\) relative to \(S_1\).

Now let \(S_2:\mathcal C\to\mathcal D\) be a shadow for which there
exist non-isomorphic objects \(X,Y\in\mathcal C\) satisfying $S_2(X)\cong S_2(Y)$.
Then the \(S_2\)-shadow does not classify \(X\) within its shadow fibre,
since that fibre also contains the non-isomorphic object \(Y\).
Consequently \(X\) cannot have lift height \(0\) relative to \(S_2\).
If the lift tower for \(S_2\) resolves this ambiguity at some finite
stage, then \(X\) has positive lift height relative to \(S_2\).

Thus \(X\) has lift height \(0\) relative to the identity shadow, while
relative to \(S_2\) its lift height is nonzero whenever it is defined
at a finite resolving stage. Therefore lift height depends on the
chosen shadow and is not, in general, an intrinsic property of the
object alone.
\end{proof}

\begin{remark}
The preceding theorem identifies commutator lift height with nilpotency
class relative to the abelianisation shadow. This should not be
interpreted as saying that lift height is merely a reformulation of
nilpotency class. Nilpotency class is an intrinsic invariant of a
group, whereas lift height is a relative notion depending on the shadow
through which hidden structure is recovered. Different shadows may
produce different lift towers and different lift heights for the same
object.
\end{remark}

\section{Relative Lift Height and Shadow Refinement}\label{sec:relative}

The commutator lift-height formula of
Theorem~\ref{thm:nilpotency-lift-height} shows that, relative to the
abelianisation shadow, lift height is governed by the lower central
series. We now develop a more general perspective in which lift height
is measured relative to arbitrary shadows and compared across different
levels of visible information.

This reveals that lift height is fundamentally a shadow-dependent notion
of classification complexity rather than an intrinsic invariant of the
underlying object alone. We further show that lift height does not
exhaust the information contained in lift packages, and that lift
packages may distinguish objects having identical classical
lower-central-series invariants.

\subsection{Relative Commutator Lift Height}

For a nilpotent group \(G\) of class \(c\), consider the family of shadow
maps
\[
S_r(G)=G/\gamma_{r+1}(G),
\qquad
1\le r\le c.
\]

These shadows retain progressively more commutator information as \(r\)
increases. The corresponding commutator lift height therefore measures
the depth of the hidden commutator structure remaining beyond the chosen
shadow.

\begin{theorem}
\label{thm:relative-height}
Let \(G\) be a nilpotent group of class \(c\), and for
\(1 \le r \le c\) consider the shadow map
\[
S_r(G)=G/\gamma_{r+1}(G).
\]
Let \(h_{S_r}^{\mathrm{com}}(G)\) denote the commutator lift height of
\(G\) relative to the shadow \(S_r\). Then
\[
h_{S_r}^{\mathrm{com}}(G)=c-r.
\]
\end{theorem}

\begin{proof}
The shadow \(S_r(G)=G/\gamma_{r+1}(G)\) retains all commutator
information of weights at most \(r\). Thus the commutator structure
hidden by the shadow begins at weight \(r+1\). The corresponding
commutator lift tower therefore consists of the successive
lower-central quotients
\[
\gamma_{r+1}(G)/\gamma_{r+2}(G),
\quad
\gamma_{r+2}(G)/\gamma_{r+3}(G),
\quad
\dots,
\quad
\gamma_c(G)/\gamma_{c+1}(G).
\]

Since \(G\) has nilpotency class \(c\), we have
\(\gamma_c(G)\neq 1\) and \(\gamma_{c+1}(G)=1\). Each of the displayed
quotients is nontrivial. Indeed, if
\(\gamma_i(G)/\gamma_{i+1}(G)\) were trivial for some
\(r+1\leq i\leq c\), then
\(\gamma_i(G)=\gamma_{i+1}(G)\), so the lower central series would
stabilise from that point onward. Since \(\gamma_{c+1}(G)=1\), this
would force \(\gamma_i(G)=1\), contradicting
\(\gamma_c(G)\neq 1\).

There are therefore exactly \(c-r\) nontrivial commutator lift layers
beyond the shadow \(S_r\). By definition of commutator lift height, the
height is equal to the number of nontrivial lift layers required to
recover the hidden commutator structure. Hence
\[
h_{S_r}^{\mathrm{com}}(G)=c-r.
\]
\end{proof}

The theorem shows that lift height decreases as more commutator
structure is made visible in the shadow. In particular, lift height
cannot be regarded as an intrinsic invariant of the underlying object
but must always be understood relative to the chosen shadow.

\begin{corollary}
\label{cor:shadow-dependent-height}
Let \(G\) be a nilpotent group of class \(c\geq 2\). Then the commutator lift height of \(G\) depends on the choice of shadow map. In particular, lift height is not, in general, 
an intrinsic invariant of the underlying object alone.
\end{corollary}

\begin{proof}
For the abelianisation shadow \(S_1(G)=G/\gamma_2(G)\), Theorem~\ref{thm:relative-height} gives \(h_{S_1}^{\mathrm{com}}(G)=c-1\). More generally, \(h_{S_r}^{\mathrm{com}}(G)=c-r\). 
Thus, for \(1\leq r,s\leq c\) with \(r\neq s\), \(h_{S_r}^{\mathrm{com}}(G)\neq h_{S_s}^{\mathrm{com}}(G)\). Hence the same group may have different lift heights relative to different shadows. 
Therefore lift height measures classification complexity relative to a chosen shadow rather than a property of the object in isolation.
\end{proof}

The preceding result shows that lift height is not an absolute
invariant of an object but depends on the amount of hidden structure
remaining beyond the chosen shadow. This raises two natural questions:
which lift heights can occur, and to what extent does lift height
determine the underlying hidden structure?

\subsection{Monotonicity Under Shadow Refinement}

The preceding examples suggest that retaining more visible information
should reduce the amount of hidden structure that must be recovered by
lifting. We formalise this intuition through a partial ordering on shadow
maps.

\begin{definition} Let $S_1:\mathcal C\to\mathcal D_1$, $S_2:\mathcal C\to\mathcal D_2$ be shadow maps. We write $S_1\preceq S_2$
if \(S_2\) retains at least as much information as \(S_1\), in the sense that there exists a functor $\pi:\mathcal D_2\to\mathcal D_1$ together with a natural isomorphism
$S_1\cong \pi\circ S_2$.
\end{definition}

Thus \(S_2\) may be regarded as a refinement of \(S_1\).

\begin{theorem}
\label{thm:monotonicity}
Let \(S_1\preceq S_2\) be shadow maps, and let
\(\{\mathcal L_n^{(1)}\}_{n\geq 0}\) and
\(\{\mathcal L_n^{(2)}\}_{n\geq 0}\) be lift towers associated with
\(S_1\) and \(S_2\), respectively. Suppose that, for every object
\(X\in\mathcal C\) and every \(n\geq 0\), the lift package
\(\mathcal L_n^{(2)}\) on the \(S_2(X)\)-shadow fibre refines the
restriction of \(\mathcal L_n^{(1)}\) to that fibre. Then, whenever
\(h_{S_1}(X)\) is finite, $h_{S_2}(X)\leq h_{S_1}(X)$. In particular, under compatible refinement of the associated lift
towers, retaining more visible information cannot increase the number
of lift layers required for reconstruction.
\end{theorem}

\begin{proof}
Since \(S_1\preceq S_2\), there exists a functor \(\pi:\mathcal D_2\to\mathcal D_1\) together with a natural isomorphism $S_1\cong \pi\circ S_2$.If \(Y\in\mathcal C\) 
satisfies \(S_2(Y)\cong S_2(X)\), then applying \(\pi\) gives
\[
S_1(Y)
\cong
\pi(S_2(Y))
\cong
\pi(S_2(X))
\cong
S_1(X).
\]
Hence $\mathcal C_{S_2(X)} \subseteq \mathcal C_{S_1(X)}$. Let \(n=h_{S_1}(X)\). By definition of lift height,
\(\mathcal L_n^{(1)}\) is complete on the \(S_1(X)\)-shadow fibre.
Its restriction to the smaller \(S_2(X)\)-shadow fibre is therefore
also complete. By hypothesis, \(\mathcal L_n^{(2)}\) refines this restricted lift
package. Proposition~\ref{prop:refinement-preserves-completeness}
therefore implies that \(\mathcal L_n^{(2)}\) is complete on the
\(S_2(X)\)-shadow fibre. Thus the \(S_2\)-lift tower reconstructs
\(X\) by stage \(n\), and hence $h_{S_2}(X)\leq n=h_{S_1}(X)$.
\end{proof}

\begin{corollary}
\label{cor:strict-shadow-monotonicity}
Let \(G\) be a nilpotent group of class \(c\geq 2\), and let
\(S_r(G)=G/\gamma_{r+1}(G)\). For integers \(1\leq r<s\leq c\),
one has
\(h_{S_r}^{\mathrm{com}}(G)>h_{S_s}^{\mathrm{com}}(G)\).
In particular, the standard lower-central-series shadows form a strictly
decreasing hierarchy of commutator lift heights
\[
c-1
>
c-2
>
\cdots
>
1
>
0.
\]
\end{corollary}

\begin{proof}
By Theorem~\ref{thm:relative-height},
\(h_{S_r}^{\mathrm{com}}(G)=c-r\) and
\(h_{S_s}^{\mathrm{com}}(G)=c-s\). Since \(r<s\leq c\),
we have \(c-r>c-s\). Hence
\(h_{S_r}^{\mathrm{com}}(G)>h_{S_s}^{\mathrm{com}}(G)\).
The displayed chain is obtained by evaluating
\(h_{S_r}^{\mathrm{com}}(G)=c-r\) for \(r=1,\dots,c\).
\end{proof}

Monotonicity shows that lift height decreases as more hidden structure
is made visible in the chosen shadow. One may therefore ask whether the
remaining lift structure is completely determined by familiar
invariants associated with that shadow. The next result shows that this
is not the case, even in the nilpotent setting.

\subsection{Lift Packages Beyond Lower-Central-Series Ranks}

A natural question is whether lift packages are determined by familiar
coarse invariants such as nilpotency class and the ranks of the
successive lower-central quotients. The following theorem shows that
they are not.

\begin{theorem}
\label{thm:beyond-classical}
There exist nilpotent Lie algebras \(L_1\) and \(L_2\) such that
\(\operatorname{cl}(L_1)=\operatorname{cl}(L_2)\), and
\[
\dim\!\bigl(\gamma_i(L_1)/\gamma_{i+1}(L_1)\bigr)
=
\dim\!\bigl(\gamma_i(L_2)/\gamma_{i+1}(L_2)\bigr)
\]
for every \(i\), yet the corresponding lift packages are not
isomorphic. Consequently, lift structure contains information not
determined by nilpotency class and the lower-central-series ranks.
\end{theorem}

\begin{proof}
Let \(F\) be a field of characteristic different from \(2\). Let \(L_1\)
have basis \(a,b,c,d,u,v\), with nonzero brackets \([a,b]=u\) and
\([c,d]=v\). Let \(L_2\) have basis \(a,b,c,d,u,v\), with nonzero
brackets \([a,b]=u\) and \([a,c]=v\). In both cases \(u\) and \(v\) are
central, and all brackets not forced by skew-symmetry are zero. Hence
both \(L_1\) and \(L_2\) are nilpotent Lie algebras of class \(2\). For both algebras, \(\gamma_2(L_i)=\langle u,v\rangle\) and
\(\gamma_3(L_i)=0\). Therefore
\(\dim(L_i/\gamma_2(L_i))=4\) and \(\dim(\gamma_2(L_i))=2\) for
\(i=1,2\). Thus \(L_1\) and \(L_2\) have the same nilpotency class and
the same lower-central-series ranks. We now compare their lift packages. For a class-two Lie algebra, the
first lift package is determined by the alternating commutator map
\(\beta_L:\Lambda^2(L/\gamma_2(L))\to\gamma_2(L)\). For \(L_1\), the
commutator map has two independent components:
\(a\wedge b\mapsto u\) and \(c\wedge d\mapsto v\). Choose a linear
functional \(\lambda:\gamma_2(L_1)\to F\) satisfying
\(\lambda(u)=\lambda(v)=1\). Then the scalar-valued alternating form
\(\lambda\circ\beta_{L_1}\) satisfies
\((\lambda\circ\beta_{L_1})(a,b)=1\) and
\((\lambda\circ\beta_{L_1})(c,d)=1\), with all other pairings between
distinct basis vectors zero except those forced by skew-symmetry. If
\(x=x_a a+x_b b+x_c c+x_d d\) lies in the radical of this form, then
pairing \(x\) with \(b,a,d,c\), respectively, gives
\(x_a=x_b=x_c=x_d=0\). Thus the form has trivial radical and hence rank
\(4\) on \(L_1/\gamma_2(L_1)\).

For \(L_2\), let \(\lambda:\gamma_2(L_2)\to F\) be any linear
functional, and write \(\alpha=\lambda(u)\) and
\(\delta=\lambda(v)\). Then the scalar-valued alternating form
\(\lambda\circ\beta_{L_2}\) has only the possibly nonzero pairings
\((\lambda\circ\beta_{L_2})(a,b)=\alpha\) and
\((\lambda\circ\beta_{L_2})(a,c)=\delta\), together with those forced
by skew-symmetry. The vector \(d\) always lies in the radical. If
\((\alpha,\delta)\neq(0,0)\), then \(\delta b-\alpha c\) gives a second
independent vector in the radical, while if \(\alpha=\delta=0\) the
form is zero. Hence every such alternating form has rank at most \(2\). An isomorphism of lift packages induces linear isomorphisms on
\(L/\gamma_2(L)\) and \(\gamma_2(L)\) intertwining the commutator maps.
Consequently, the family of scalar-valued alternating forms obtained by
composing \(\beta_L\) with linear functionals on \(\gamma_2(L)\) is
preserved up to change of basis. In particular, the maximal rank among
such forms is an invariant of the lift package. For \(L_1\) this
maximal rank is \(4\), while for \(L_2\) it is at most \(2\). Hence the
lift packages of \(L_1\) and \(L_2\) are not isomorphic. Therefore the
lift package is not determined by nilpotency class and
lower-central-series ranks alone.
\end{proof}

\medskip

The results of this section show that lift height is relative to a chosen
shadow, behaves monotonically under shadow refinement, and does not
exhaust the information contained in lift packages. Thus lift theory
captures both the depth and the organisation of the hidden structure
omitted by a shadow.

\section{Canonical Lift Towers and the Depth of Classification Failure}
\label{sec:canonical-towers}

A fixed lift tower induces a descending chain of equivalence relations
on each shadow fibre, recording the successive distinctions made by
its lift layers. Such a tower may contain redundant stages which add
no new classification information. Canonical reduction removes these
repetitions, while essential lift height records the resulting strict
refinement depth.

This separates the formal length of a lift tower from the depth at
which genuinely new distinctions occur. The second part of the section
asks whether such new classification information can occur arbitrarily
deep in a commutator tower. We show that it can, even over a fixed
four-dimensional abelianisation shadow and even when the underlying
lower-central graded vector space is fixed.

\subsection{Canonical Reduction and Essential Lift Height}

Rather than working only with the individual lift operations, it is
useful to record the sequence of equivalence relations that they induce
on a shadow fibre. This refinement chain makes visible which stages
produce genuinely new distinctions and which merely repeat information
already present at an earlier stage.

\begin{definition}
\label{def:refinement-chain}
The \emph{refinement chain} of $\mathcal T$ is the ordered sequence
\[
R(\mathcal T)
=
\bigl(
\sim_0^{\mathcal T},
\sim_1^{\mathcal T},
\ldots,
\sim_n^{\mathcal T}
\bigr).
\]
\end{definition}

\begin{definition}
The layer $\mathcal L_i$ is called \emph{redundant} if $\sim_{i-1}^{\mathcal T} = \sim_i^{\mathcal T}$.
Equivalently, the first $i$ layers distinguish exactly the same pairs
of objects as the first $i-1$ layers.
\end{definition}

\begin{definition}
\label{def:reduced-lift-tower}
A lift tower is called \emph{reduced} if $\sim_{i-1}^{\mathcal T} \supsetneq \sim_i^{\mathcal T}$ for every $1\leq i\leq n$.
\end{definition}

\begin{definition}
\label{def:refinement-equivalent}
Two finite lift towers $\mathcal T$ and $\mathcal T'$ on the same
shadow fibre are called \emph{refinement-equivalent} if the strict
chains obtained from $R(\mathcal T)$ and $R(\mathcal T')$ by deleting
repeated consecutive relations are identical.
\end{definition}

\begin{remark}
Refinement-equivalence does not require the individual lift operations
in the two towers to be isomorphic. It records only whether the towers
induce the same essential sequence of refinements on the shadow fibre.
Consequently, refinement-equivalence is weaker than an isomorphism of
lift packages.
\end{remark}

\begin{theorem}
\label{thm:canonical-reduction}
Let $\mathcal T$ be a finite lift tower with refinement chain
$R(\mathcal T)=\bigl(\sim_0^{\mathcal T},\sim_1^{\mathcal T},\ldots,\sim_n^{\mathcal T}\bigr)$.
Then there exists a unique strict chain
$\operatorname{Red}(\mathcal T)=\bigl(\rho_0,\rho_1,\ldots,\rho_r\bigr)$
obtained from $R(\mathcal T)$ by deleting repeated consecutive
equivalence relations. Moreover,
$\rho_0\supsetneq\rho_1\supsetneq\cdots\supsetneq\rho_r$.
\end{theorem}

\begin{proof}
Since $R(\mathcal T)$ is a refinement chain, we have
$\sim_0^{\mathcal T}\supseteq\sim_1^{\mathcal T}\supseteq\cdots\supseteq\sim_n^{\mathcal T}$.

We first show that repeated equivalence relations occur only in
consecutive blocks. Suppose that
$\sim_i^{\mathcal T}=\sim_k^{\mathcal T}$
for some indices $0\leq i<k\leq n$. Then, for every
$j$ satisfying $i\leq j\leq k$,
$\sim_i^{\mathcal T}\supseteq\sim_j^{\mathcal T}\supseteq\sim_k^{\mathcal T}=\sim_i^{\mathcal T}$.
Hence both inclusions are equalities, so
$\sim_i^{\mathcal T}=\sim_j^{\mathcal T}=\sim_k^{\mathcal T}$.
Therefore every equivalence relation appearing in
$R(\mathcal T)$ occurs in a unique maximal consecutive block of equal
terms. Replace each maximal constant block by a single representative. This
produces a chain
$\operatorname{Red}(\mathcal T)=(\rho_0,\rho_1,\ldots,\rho_r)$.
Since consecutive blocks correspond to distinct equivalence relations,
the terms $\rho_0,\ldots,\rho_r$ are pairwise distinct. Furthermore,
because the original chain is descending,
$\rho_j\supseteq\rho_{j+1}$ for $0\leq j<r$.
As $\rho_j\neq\rho_{j+1}$, each inclusion is strict. Hence
$\rho_0\supsetneq\rho_1\supsetneq\cdots\supsetneq\rho_r$.
This proves existence.

To prove uniqueness, let
$(\sigma_0,\sigma_1,\ldots,\sigma_s)$
be any strict chain obtained from $R(\mathcal T)$ by deleting repeated
consecutive terms. Since only repeated terms are removed, each maximal
constant block contributes at least one representative to the resulting
chain. On the other hand, because the chain is strict, no block can
contribute more than one representative. Consequently each maximal
constant block contributes exactly one representative, and these appear
in the same order as in the original refinement chain. Since every representative of a given block is the same equivalence
relation, the resulting strict chain is uniquely determined. Therefore
$(\sigma_0,\sigma_1,\ldots,\sigma_s)=(\rho_0,\rho_1,\ldots,\rho_r)$,
and hence
$\operatorname{Red}(\mathcal T)=(\rho_0,\rho_1,\ldots,\rho_r)$
is unique.
\end{proof}

\begin{remark}
The canonical object produced by
Theorem~\ref{thm:canonical-reduction} is the reduced refinement chain
of a fixed lift tower, not a uniquely determined reduced lift tower.
Thus the theorem does not assert that arbitrary reduced lift towers
over the same shadow fibre are refinement-equivalent. Different lift
towers may encode genuinely different sequences of strict refinements.
\end{remark}

\begin{definition}
\label{def:essential-height}
Let $\mathcal T$ be a finite lift tower and write $\operatorname{Red}(\mathcal T) = (\rho_0,\rho_1,\ldots,\rho_r)$.
The \emph{essential height} of $\mathcal T$ is $h_{\mathrm{ess}}(\mathcal T)=r$. For an object $X$, one may similarly define its essential height to be
the first index at which the reduced refinement chain determines the
isomorphism class of $X$ within its shadow fibre.
\end{definition}

\begin{remark}
\label{rem:essential-height-bound}
Essential lift height is unchanged by the insertion or deletion of
redundant stages, since these operations do not change
$\operatorname{Red}(\mathcal T)$. It therefore depends only on the
refinement-equivalence class of the fixed tower, but need not be
independent of the choice of lift tower.

For a nilpotent group \(G\) of class \(c\) and its \(r\)-step
nilpotent shadow \(S_r(G)\), where \(1\le r<c\), this gives
\(h_{\mathrm{ess}}(\mathcal T^{\mathrm{com}}_{S_r(G)})
\le h_{S_r}^{\mathrm{com}}(G)=c-r\). In particular,
\(h_{\mathrm{ess}}(\mathcal T^{\mathrm{com}}_{S_1(G)})
\le h_{\mathrm{com}}(G)=c-1\).

Equality
\(h_{\mathrm{ess}}(\mathcal T^{\mathrm{com}}_{S_r(G)})=c-r\)
holds if and only if the commutator lift tower above \(S_r(G)\) is
reduced, equivalently if and only if every commutator lift stage gives
a strict refinement of the induced equivalence relation.
\end{remark}

Essential height records how many stages of a tower make genuinely new
distinctions. The complementary question is what can be recovered
before such a distinction appears. If two objects remain equivalent
after a given truncation, then no invariant depending only on that
truncation can separate them.

\begin{proposition}
\label{prop:truncation-obstruction}
Let $\mathcal T=(\mathcal L_1,\ldots,\mathcal L_n)$ be a lift tower on
a shadow fibre $\mathcal C_A$, and let $0\leq k\leq n$. Write
$\mathcal T_{\leq k}=(\mathcal L_1,\ldots,\mathcal L_k)$ for its
$k$-th truncation. Suppose that $X\sim_k^{\mathcal T}Y$ for objects
$X,Y\in\mathcal C_A$. Let $F:\mathcal C_A\longrightarrow\mathcal E$ be an invariant which
factors through $\mathcal T_{\leq k}$. Thus there exists an invariant
$\overline F:\operatorname{Im}(\mathcal T_{\leq k})\longrightarrow
\mathcal E$ such that
$F=\overline F\circ\mathcal T_{\leq k}$. Then $F(X)\cong F(Y)$.
Consequently, no invariant depending only upon the first $k$ lift
layers can distinguish $X$ from $Y$.
\end{proposition}

\begin{proof}
By the definition of the equivalence relation $\sim_k^{\mathcal T}$,
$X\sim_k^{\mathcal T}Y$ if and only if
$\mathcal T_{\leq k}(X)\cong\mathcal T_{\leq k}(Y)$. Since $F$ factors
through $\mathcal T_{\leq k}$,
$F=\overline F\circ\mathcal T_{\leq k}$. Therefore
$F(X)=\overline F\!\left(\mathcal T_{\leq k}(X)\right)
\cong\overline F\!\left(\mathcal T_{\leq k}(Y)\right)=F(Y)$,
where the middle isomorphism follows because $\overline F$ preserves
isomorphisms. Hence $F$ assigns isomorphic values to $X$ and $Y$, so it
cannot distinguish them.
\end{proof}

Proposition~\ref{prop:truncation-obstruction} is formal: it explains
what follows once two objects agree through a prescribed truncation.
The substantive question is whether pairs satisfying
$X\sim_k^{\mathcal T}Y$ but
$X\not\sim_{k+1}^{\mathcal T}Y$ actually occur. We first exhibit a
class-three example in which both successive commutator lift stages are
essential, and then show that the same phenomenon can occur at
arbitrarily large depth.

\subsection{Essential Commutator Layers}

The first example shows that essential height is not merely a device
for removing artificial repetitions: successive commutator layers can
each contribute genuinely new classification information. A
class-three tower is the smallest setting in which two successive
strict refinements can occur.

\begin{theorem}
\label{thm:arbitrarily-deep-essential-layers}
Let \(m\geq 3\), and let \(F\) be a field of characteristic different
from \(2\). Then there exist nilpotent Lie algebras \(L_m\) and \(M_m\)
of nilpotency class exactly \(m\), with
\(L_m/\gamma_2(L_m)\cong M_m/\gamma_2(M_m)\cong F^4\), such that their
commutator lift packages agree through degree \(m-1\) but are not
isomorphic after the degree-\(m\) operation is included. Equivalently, in the induced commutator lift tower,
\(L_m\sim_{m-2}M_m\) and \(L_m\not\sim_{m-1}M_m\). Thus the final
commutator lift layer can be essential at arbitrarily large nilpotency
class.
\end{theorem}

\begin{proof}
Let \(V_1=\langle a,b,c,d\rangle\),
\(V_i=\langle x_i,y_i\rangle\) for \(2\leq i\leq m-1\), and
\(V_m=\langle z\rangle\), and let
\[V=V_1\oplus V_2\oplus\cdots\oplus V_m.\]

We define two Lie algebra structures \(L_m\) and \(M_m\) on \(V\).
In both algebras set \([a,b]=x_2\) and \([c,d]=y_2\), and, for
\(2\leq i\leq m-2\), set \([a,x_i]=x_{i+1}\) and
\([c,y_i]=y_{i+1}\). For \(L_m\), define the final nonzero bracket by
\([a,x_{m-1}]=z\). For \(M_m\), define
\([a,x_{m-1}]=z\) and \([c,y_{m-1}]=z\). All remaining brackets
between basis elements are zero, apart from those forced by
skew-symmetry. We first verify the Jacobi identity. The nonzero brackets split into
two chains,
\([a,b]=x_2,\ [a,x_2]=x_3,\ldots,[a,x_{m-1}]=z\),
and
\([c,d]=y_2,\ [c,y_2]=y_3,\ldots\).
In \(M_m\), the second chain also has the final bracket
\([c,y_{m-1}]=z\). There are no nonzero cross-brackets between the two
chains, and \(z\) is central.

Within the first chain, every nonzero bracket involves \(a\), while
within the second chain every nonzero bracket involves \(c\). Hence,
for a triple of basis elements contained in one chain, every Jacobi
sum either vanishes termwise or reduces to a bracket of two elements
among the \(x_i\), respectively among the \(y_i\), and is therefore
zero. A triple involving elements from both chains also has zero
Jacobi sum because all cross-brackets vanish. Thus both \(L_m\) and
\(M_m\) satisfy the Jacobi identity. Their lower central series may be read directly from the defining
brackets. For \(2\leq i\leq m-1\),
\(\gamma_i(L_m)=
\langle x_i,\ldots,x_{m-1},y_i,\ldots,y_{m-1},z\rangle\),
and the same equality holds for \(M_m\). At the final two stages,
\(\gamma_m(L_m)=\gamma_m(M_m)=\langle z\rangle\) and
\(\gamma_{m+1}(L_m)=\gamma_{m+1}(M_m)=0\). Therefore both algebras
have nilpotency class exactly \(m\).

Their abelianisations are naturally identified with
\(A=\langle\bar a,\bar b,\bar c,\bar d\rangle\cong F^4\). Moreover,
for \(2\leq i\leq m-1\),
\(K_i=\gamma_i/\gamma_{i+1}=\langle\bar x_i,\bar y_i\rangle\), while
\(K_m=\gamma_m=\langle z\rangle\). By construction, \(L_m\) and \(M_m\) have identical graded commutator
operations through degree \(m-1\). Indeed, their common operations are
determined by
\([\bar a,\bar b]=\bar x_2\),
\([\bar c,\bar d]=\bar y_2\),
\([\bar a,\bar x_i]=\bar x_{i+1}\), and
\([\bar c,\bar y_i]=\bar y_{i+1}\)
for every \(2\leq i\leq m-2\). Hence their commutator lift packages
agree through degree \(m-1\), and therefore
\(L_m\sim_{m-2}M_m\).

It remains to distinguish the final lift layer. Consider the
degree-\(m\) operation
\(\beta_m:A\otimes K_{m-1}\longrightarrow K_m\)
and associate to it the linear map
\(\Phi_m:A\longrightarrow\operatorname{Hom}(K_{m-1},K_m)\), where
\(\Phi_m(v)(w)=\beta_m(v,w)\). For \(L_m\), the only nonzero bracket contributing to \(\beta_m\) is
\([a,x_{m-1}]=z\). Consequently,
\(\operatorname{rank}\Phi_m^{L_m}=1\). For \(M_m\), the final brackets are
\([a,x_{m-1}]=z\) and \([c,y_{m-1}]=z\). Thus
\(\Phi_m^{M_m}(\bar a)\) and \(\Phi_m^{M_m}(\bar c)\) are the two
independent coordinate functionals on
\(K_{m-1}=\langle\bar x_{m-1},\bar y_{m-1}\rangle\). Therefore
\(\operatorname{rank}\Phi_m^{M_m}=2\). The rank of \(\Phi_m\) is invariant under isomorphism of the
degree-\(m\) commutator lift package. Indeed, if
\(f:A\longrightarrow A'\),
\(g:K_{m-1}\longrightarrow K'_{m-1}\), and
\(h:K_m\longrightarrow K'_m\) are isomorphisms intertwining the
degree-\(m\) operations, then
\(\Phi'_m(f(v))=h\circ\Phi_m(v)\circ g^{-1}\).
Composition with linear isomorphisms preserves rank. Hence isomorphic
degree-\(m\) lift packages have the same value of
\(\operatorname{rank}\Phi_m\). Since \(\operatorname{rank}\Phi_m^{L_m}=1\) and
\(\operatorname{rank}\Phi_m^{M_m}=2\), the degree-\(m\) commutator
lift packages of \(L_m\) and \(M_m\) are not isomorphic. Therefore
\(L_m\not\sim_{m-1}M_m\). We have thus constructed, for every \(m\geq3\), a pair of class-\(m\)
nilpotent Lie algebras satisfying \(L_m\sim_{m-2}M_m\) and
\(L_m\not\sim_{m-1}M_m\). Hence the final commutator lift layer can be
essential at arbitrarily large nilpotency class.
\end{proof}

\begin{theorem} 
\label{thm:class-three-essential-height} 
Let \(p>3\) be a prime. There exists a family \(\mathcal F\) of finite 
\(p\)-groups of nilpotency class three, all lying over the fixed 
abelianisation shadow \(A=(\mathbf F_p)^4\), for which the commutator 
lift tower is reduced. Consequently, 
\(h_{\mathrm{ess}}(\mathcal T^{\mathrm{com}}_{\mathcal F})=2\). More precisely, there exist \(G_P,G_Q,G_S\in\mathcal F\) such that 
\(G_P\sim^{\mathcal T^{\mathrm{com}}}_0 G_Q\) and 
\(G_P\not\sim^{\mathcal T^{\mathrm{com}}}_1 G_Q\), while 
\(G_Q\sim^{\mathcal T^{\mathrm{com}}}_1 G_S\) and 
\(G_Q\not\sim^{\mathcal T^{\mathrm{com}}}_2 G_S\). Thus both stages 
of the class-three commutator lift tower give strict refinements of the 
induced equivalence relation. Moreover, no invariant factoring through the first commutator lift 
package is a complete classifier on \(\mathcal F\). 
\end{theorem} 
 
\begin{proof} 
Let \(F=\mathbf F_p\). We first construct three nilpotent Lie algebras 
\(P,Q,S\) over \(F\), each with basis \(a,b,c,d,u,v,z\). All brackets 
not specified below are zero, apart from those forced by 
skew-symmetry. Define \(P\) by \([a,b]=u\), \([c,d]=v\), and \([a,u]=z\). Define 
\(Q\) by \([a,b]=u\), \([a,c]=v\), and \([a,u]=z\). Finally, define 
\(S\) by \([a,b]=u\), \([a,c]=v\), \([b,u]=z\), and \([c,v]=z\). 

These brackets satisfy the Jacobi identity. Indeed, in each case, if 
\(V_1=\langle a,b,c,d\rangle\), \(V_2=\langle u,v\rangle\), and 
\(V_3=\langle z\rangle\), then 
\([V_1,V_1]\subseteq V_2\), \([V_1,V_2]\subseteq V_3\), 
\([V_2,V_2]=0\), and \([V_1,V_3]=0\). If a Jacobi triple contains an 
element of \(V_3\), every term vanishes because \(V_3\) is central. If 
it contains at least two elements of \(V_2\), every term also vanishes 
since \([V_2,V_2]=0\) and \([V_2,V_3]=0\). If it contains exactly one 
element of \(V_2\) and two elements of \(V_1\), the inner bracket of 
the two \(V_1\)-elements lies in \(V_2\), while an inner bracket 
involving the \(V_2\)-element lies in \(V_3\), so again every Jacobi 
term vanishes. It remains only to check triples of elements of \(V_1\). 

For \(Q\) and \(S\), the only potentially nonzero Jacobi relation among 
distinct noncentral degree-one generators is 
\([a,[b,c]]+[b,[c,a]]+[c,[a,b]] 
=[a,0]+[b,-v]+[c,u]=0\). 
The generator \(d\) is central in both \(Q\) and \(S\), so every 
remaining degree-one Jacobi relation is zero. For \(P\), the four triples of distinct degree-one basis elements give 
\([a,[b,c]]+[b,[c,a]]+[c,[a,b]]=[c,u]=0\), 
\([a,[b,d]]+[b,[d,a]]+[d,[a,b]]=[d,u]=0\), 
\([a,[c,d]]+[c,[d,a]]+[d,[a,c]]=[a,v]=0\), and 
\([b,[c,d]]+[c,[d,b]]+[d,[b,c]]=[b,v]=0\). 
Thus \(P,Q,S\) are Lie algebras. For each \(L\in\{P,Q,S\}\), 
\(\gamma_2(L)=\langle u,v,z\rangle\), 
\(\gamma_3(L)=\langle z\rangle\), and \(\gamma_4(L)=0\). 
Thus all three Lie algebras have nilpotency class exactly three, and 
\(L/\gamma_2(L)\cong 
\langle\bar a,\bar b,\bar c,\bar d\rangle\cong F^4\). 
 
Since \(p>3\), the nilpotency class is strictly less than \(p\). 
The Lazard correspondence therefore associates to \(P,Q,S\) finite 
\(p\)-groups \(G_P,G_Q,G_S\) of nilpotency class three. By 
Theorem~\ref{thm:lazard}, their lower-central quotients and commutator 
lift operations correspond to the graded Lie operations of 
\(P,Q,S\), respectively. Fix \(A=(\mathbf F_p)^4\) and identify 
\(G_L/\gamma_2(G_L)\cong A\) for \(L\in\{P,Q,S\}\) using the images 
of \(a,b,c,d\). Thus \(G_P,G_Q,G_S\) lie over the same fixed 
abelianisation shadow. 

We first show that the first commutator lift stage gives a strict 
refinement. On the Lie side, the first commutator map for \(P\), 
\(\beta_{2,P}: 
\Lambda^2(P/\gamma_2(P))\longrightarrow 
\gamma_2(P)/\gamma_3(P)\), is determined by 
\(\beta_{2,P}(\bar a,\bar b)=\bar u\) and 
\(\beta_{2,P}(\bar c,\bar d)=\bar v\). For \(Q\), 
\(\beta_{2,Q}(\bar a,\bar b)=\bar u\) and 
\(\beta_{2,Q}(\bar a,\bar c)=\bar v\). These first lift packages are not isomorphic. Let 
\(K_2=\langle\bar u,\bar v\rangle\). For a linear functional 
\(\lambda:K_2\longrightarrow F\), consider the scalar alternating form 
\(\lambda\circ\beta_2\) on the abelianisation. 
 
For \(P\), choose \(\lambda\) with 
\(\lambda(\bar u)=\lambda(\bar v)=1\). Then the resulting alternating 
form has nonzero pairings \((\bar a,\bar b)=1\) and 
\((\bar c,\bar d)=1\), and zero radical. Hence it has rank \(4\). For \(Q\), write \(\alpha=\lambda(\bar u)\) and 
\(\delta=\lambda(\bar v)\). The only possibly nonzero pairings are 
\((\bar a,\bar b)=\alpha\) and 
\((\bar a,\bar c)=\delta\). The vector \(\bar d\) always lies in the 
radical. If \((\alpha,\delta)\neq(0,0)\), then 
\(\delta\bar b-\alpha\bar c\) is a second linearly independent radical 
vector, while if \(\alpha=\delta=0\), the form is zero. Thus every 
scalar specialization of \(\beta_{2,Q}\) has rank at most \(2\). The maximal rank of a scalar specialization is invariant under 
isomorphism of the first commutator lift package. Therefore the first 
lift packages associated with \(P\) and \(Q\) are not isomorphic. 
By Theorem~\ref{thm:lazard}, the same is true of the first commutator 
lift packages of \(G_P\) and \(G_Q\). Since their abelianisation 
shadows are both identified with \(A\), 
\(G_P\sim^{\mathcal T^{\mathrm{com}}}_0G_Q\) and 
\(G_P\not\sim^{\mathcal T^{\mathrm{com}}}_1G_Q\). 
Hence the first refinement is strict. 

We next show that the second commutator lift stage is strict. The first 
commutator lift packages associated with \(Q\) and \(S\) are identical: 
in both cases \(\beta_2(\bar a,\bar b)=\bar u\) and 
\(\beta_2(\bar a,\bar c)=\bar v\), with all other basic values zero. 
The Lazard correspondence therefore gives 
\(G_Q\sim^{\mathcal T^{\mathrm{com}}}_1G_S\).  
 
It remains to distinguish their second lift layers. On the Lie side 
these are described by 
\(\beta_3:A\otimes K_2\longrightarrow K_3\), where 
\(A=L/\gamma_2(L)\), 
\(K_2=\gamma_2(L)/\gamma_3(L)\), and 
\(K_3=\gamma_3(L)=\langle z\rangle\). Associate to \(\beta_3\) the 
linear map 
\(\Phi_{\beta_3}:A\longrightarrow\operatorname{Hom}(K_2,K_3)\), where 
\(\Phi_{\beta_3}(x)(y)=\beta_3(x,y)\). For \(Q\), the only nonzero degree-three bracket is \([a,u]=z\). 
Consequently, \(\operatorname{rank}\Phi_{\beta_{3,Q}}=1\). For \(S\), \([b,u]=z\) and \([c,v]=z\). Thus 
\(\Phi_{\beta_{3,S}}(\bar b)\) and 
\(\Phi_{\beta_{3,S}}(\bar c)\) are linearly independent elements of 
\(\operatorname{Hom}(K_2,K_3)\), and hence 
\(\operatorname{rank}\Phi_{\beta_{3,S}}=2\). 

This rank is invariant under isomorphism of the full commutator lift 
package. Indeed, suppose \(f:A\longrightarrow A'\), 
\(g:K_2\longrightarrow K'_2\), and 
\(h:K_3\longrightarrow K'_3\) are isomorphisms intertwining the 
degree-three operations. Then 
\(\Phi_{\beta'_3}(f(x)) 
=h\circ\Phi_{\beta_3}(x)\circ g^{-1}\). 
The map 
\(T\longmapsto h\circ T\circ g^{-1}\) is a linear isomorphism from 
\(\operatorname{Hom}(K_2,K_3)\) to 
\(\operatorname{Hom}(K'_2,K'_3)\). Since \(f\) is also a linear 
isomorphism, the images of \(\Phi_{\beta_3}\) and 
\(\Phi_{\beta'_3}\) have the same dimension. Hence 
\(\operatorname{rank}\Phi_{\beta'_3} 
=\operatorname{rank}\Phi_{\beta_3}\). 
 
But \(\operatorname{rank}\Phi_{\beta_{3,Q}}=1\) and 
\(\operatorname{rank}\Phi_{\beta_{3,S}}=2\). Hence the full 
class-three lift packages associated with \(Q\) and \(S\) are not 
isomorphic. By Theorem~\ref{thm:lazard}, the full commutator lift 
packages of \(G_Q\) and \(G_S\) are therefore not isomorphic. 
Consequently, 
\(G_Q\not\sim^{\mathcal T^{\mathrm{com}}}_2G_S\). 

Set \(\mathcal F=\{G_P,G_Q,G_S\}\). We have exhibited 
\(G_P\sim^{\mathcal T^{\mathrm{com}}}_0G_Q\) and 
\(G_P\not\sim^{\mathcal T^{\mathrm{com}}}_1G_Q\), and 
\(G_Q\sim^{\mathcal T^{\mathrm{com}}}_1G_S\) and 
\(G_Q\not\sim^{\mathcal T^{\mathrm{com}}}_2G_S\). Therefore 
\({\sim}^{\mathcal T^{\mathrm{com}}}_0 
\supsetneq{\sim}^{\mathcal T^{\mathrm{com}}}_1 
\supsetneq{\sim}^{\mathcal T^{\mathrm{com}}}_2\) 
on \(\mathcal F\). Thus the two-stage commutator lift tower on 
\(\mathcal F\) is already reduced. Its canonical reduction has two 
strict stages, so 
\(h_{\mathrm{ess}}(\mathcal T^{\mathrm{com}}_{\mathcal F})=2\). Finally, 
\(G_Q\sim^{\mathcal T^{\mathrm{com}}}_1G_S\) but 
\(G_Q\not\sim^{\mathcal T^{\mathrm{com}}}_2G_S\). 
By Proposition~\ref{prop:truncation-obstruction}, every invariant 
factoring through the first commutator lift package takes isomorphic 
values on \(G_Q\) and \(G_S\). Since their full commutator lift 
packages are nonisomorphic, no invariant factoring through the first 
commutator lift package can be a complete classifier on 
\(\mathcal F\). 
\end{proof}

The class-three construction shows that consecutive commutator lift
stages can both be essential. There is, however, no bounded-depth
version of this phenomenon. By extending the same basic mechanism
through longer lower-central chains, one can postpone the first new
distinction to an arbitrarily late commutator layer.

\subsection{Unbounded Classification Depth over a Fixed Shadow}

The preceding construction varies the nilpotency class in order to
place an essential layer arbitrarily deep in the commutator tower.
Crucially, it does not require the visible shadow to grow: every member
of the family has abelianisation $F^4$. Thus the failure of bounded
classification depth already occurs inside a single fixed shadow
fibre.

\begin{theorem}
\label{thm:unbounded-classification-depth-fixed-shadow}
Let \(F\) be a field of characteristic different from \(2\), and let
\(A=F^4\). For every integer \(N\geq 1\), there exist nilpotent Lie
algebras \(L\) and \(M\) in the abelianisation shadow fibre over \(A\)
whose first \(N\) commutator lift layers are isomorphic, but whose
commutator lift packages cease to be isomorphic at the next layer. Equivalently, for every \(N\geq1\), there exist nilpotent Lie algebras
\(L\) and \(M\) satisfying
\(L/\gamma_2(L)\cong M/\gamma_2(M)\cong A\) and
\(L\sim_N M\), but \(L\not\sim_{N+1}M\).

In particular, classification failure over the fixed shadow \(A=F^4\)
can persist through arbitrarily long finite commutator lift towers.
\end{theorem}

\begin{proof}
Fix \(N\geq1\), and apply
Theorem~\ref{thm:arbitrarily-deep-essential-layers} with
\(m=N+2\). The theorem gives nilpotent Lie algebras \(L_m\) and
\(M_m\) of class exactly \(m\) satisfying
\(L_m/\gamma_2(L_m)\cong M_m/\gamma_2(M_m)\cong F^4\), together with
\(L_m\sim_{m-2}M_m\) and \(L_m\not\sim_{m-1}M_m\). Since \(m=N+2\), we have \(m-2=N\) and \(m-1=N+1\). Hence
\(L_m\sim_NM_m\), while \(L_m\not\sim_{N+1}M_m\). After identifying
their abelianisations with \(A=F^4\), both algebras therefore lie in
the same fixed abelianisation shadow fibre, agree through the first
\(N\) commutator lift layers, and are distinguished by the next
layer. Since \(N\) was arbitrary, there is no finite bound on the number of
successive commutator lift layers through which classification failure
may persist over the fixed shadow \(A=F^4\).
\end{proof}

\begin{corollary}
\label{cor:no-fixed-shadow-depth-bound}
Let \(F\) be a field of characteristic different from \(2\), and let
\(A=F^4\). There is no integer \(B\) such that, for every nilpotent
Lie algebra in the abelianisation shadow fibre over \(A\), the first
\(B\) commutator lift layers determine its full commutator lift
package up to isomorphism.
\end{corollary}

\begin{proof}
Suppose that such an integer \(B\) existed. By
Theorem~\ref{thm:unbounded-classification-depth-fixed-shadow}, with
\(N=B\), there would exist nilpotent Lie algebras \(L\) and \(M\) in
the shadow fibre over \(A\) satisfying \(L\sim_BM\) but
\(L\not\sim_{B+1}M\). The relation \(L\sim_BM\) says that the first \(B\) commutator lift
layers of \(L\) and \(M\) are isomorphic. By the assumed property of
\(B\), those first \(B\) layers would determine the full commutator
lift package, so the full packages of \(L\) and \(M\) would have to be
isomorphic. In particular their truncations through layer \(B+1\)
would be isomorphic, giving \(L\sim_{B+1}M\). This contradicts
\(L\not\sim_{B+1}M\). Therefore no such finite bound exists.
\end{proof}

\begin{corollary}
\label{cor:no-shadow-dimension-bound}
There is no function \(B:\mathbb N\to\mathbb N\) with the property
that, for every field \(F\) of characteristic different from \(2\)
and every nilpotent Lie algebra \(L\) over \(F\), the first
\(B(\dim(L/\gamma_2(L)))\) commutator lift layers determine the full
commutator lift package of \(L\) within its abelianisation shadow
fibre.
\end{corollary}

\begin{proof}
Suppose that such a function \(B\) existed. Fix a field \(F\) of
characteristic different from \(2\), and put \(N=B(4)\). By
Theorem~\ref{thm:unbounded-classification-depth-fixed-shadow}, there
exist nilpotent Lie algebras \(L\) and \(M\) satisfying
\(L/\gamma_2(L)\cong M/\gamma_2(M)\cong F^4\) and \(L\sim_NM\), but
\(L\not\sim_{N+1}M\). Both abelianisations have dimension \(4\), so the assumed property of
\(B\) says that their first \(B(4)=N\) commutator lift layers determine
their full commutator lift packages within this shadow fibre. Since
\(L\sim_NM\), those first \(N\) layers are isomorphic, and hence the
full commutator lift packages would have to be isomorphic. Their
truncations through layer \(N+1\) would then also be isomorphic,
contradicting \(L\not\sim_{N+1}M\). Thus no such function \(B\) exists. In fact, the preceding theorem
shows more: the obstruction already occurs over the single fixed
four-dimensional shadow \(F^4\).
\end{proof}

These results rule out bounds depending only on the visible
abelianisation shadow. One may ask whether substantially more
lower-central information restores such a bound. The same examples
show that it does not: the underlying graded vector spaces can be
fixed degree by degree while new commutator structure first appears
arbitrarily far along the grading.

\begin{theorem}
\label{thm:graded-vector-space-arbitrary-depth}
Let \(m\geq3\), and let \(F\) be a field of characteristic different
from \(2\). Then there exist nilpotent Lie algebras \(L_m\) and \(M_m\)
of class exactly \(m\) such that their associated lower-central graded
objects are isomorphic as graded vector spaces, while they are not
isomorphic as graded Lie algebras. Moreover, their graded Lie brackets
agree through degree \(m-1\), and the failure of isomorphism is first
detected by the degree-\(m\) commutator operation.
\end{theorem}

\begin{proof}
Take the Lie algebras \(L_m\) and \(M_m\) constructed in
Theorem~\ref{thm:arbitrarily-deep-essential-layers}. In both cases the
degree-one lower-central quotient is naturally identified with
\(F^4\). For every \(2\leq i\leq m-1\), the degree-\(i\) quotient
\(K_i=\gamma_i/\gamma_{i+1}\) has basis
\(\{\bar x_i,\bar y_i\}\), and the final quotient \(K_m=\gamma_m\)
has basis \(\{z\}\). Hence the associated graded objects of \(L_m\)
and \(M_m\) have the same homogeneous components: degree \(1\) has
dimension \(4\), every degree \(2,\ldots,m-1\) has dimension \(2\),
and degree \(m\) has dimension \(1\). Using the indicated bases gives
an isomorphism
\(\operatorname{gr}(L_m)\cong\operatorname{gr}(M_m)\) as graded
vector spaces. The defining brackets of the two algebras are identical through degree
\(m-1\). Thus the induced graded Lie brackets agree through that
degree. The difference occurs in the final operation
\(\beta_m:K_1\otimes K_{m-1}\to K_m\). For \(L_m\), the only nonzero
final bracket is \([a,x_{m-1}]=z\), whereas for \(M_m\) the final
brackets are \([a,x_{m-1}]=z\) and \([c,y_{m-1}]=z\).

As in Theorem~\ref{thm:arbitrarily-deep-essential-layers}, associate
to the final bracket the linear map
\(\Phi_m:K_1\to\operatorname{Hom}(K_{m-1},K_m)\), defined by
\(\Phi_m(v)(w)=\beta_m(v,w)\). For \(L_m\), the image of \(\Phi_m\)
is one-dimensional, since only the class of \(a\) acts nontrivially
on \(K_{m-1}\). Thus
\(\operatorname{rank}\Phi_m^{L_m}=1\). For \(M_m\), the classes of
\(a\) and \(c\) induce the two independent coordinate functionals on
\(K_{m-1}=\langle\bar x_{m-1},\bar y_{m-1}\rangle\), and therefore
\(\operatorname{rank}\Phi_m^{M_m}=2\). Suppose, for contradiction, that there were a graded Lie algebra
isomorphism
\(\psi:\operatorname{gr}(L_m)\to\operatorname{gr}(M_m)\). Since
\(\psi\) preserves degree, it induces linear isomorphisms
\(f:K_1^{L_m}\to K_1^{M_m}\),
\(g:K_{m-1}^{L_m}\to K_{m-1}^{M_m}\), and
\(h:K_m^{L_m}\to K_m^{M_m}\). Since \(\psi\) preserves the graded Lie
bracket, these maps intertwine the final operations, and consequently
\(\Phi_m^{M_m}(f(v))=h\circ\Phi_m^{L_m}(v)\circ g^{-1}\) for every
\(v\in K_1^{L_m}\).

Precomposition and postcomposition by linear isomorphisms preserve
the dimension of the image of \(\Phi_m\). Hence a graded Lie algebra
isomorphism would imply
\(\operatorname{rank}\Phi_m^{L_m}
=\operatorname{rank}\Phi_m^{M_m}\). This is impossible, since these
ranks are \(1\) and \(2\), respectively. Therefore
\(\operatorname{gr}(L_m)\) and
\(\operatorname{gr}(M_m)\) are isomorphic as graded vector spaces but
not as graded Lie algebras. Since all graded brackets through degree
\(m-1\) agree, the distinction is first detected by the degree-\(m\)
commutator operation. As \(m\) is arbitrary, this separation can occur
at arbitrarily large degree.
\end{proof}

\begin{corollary}
\label{cor:graded-dimensions-do-not-control-depth}
The dimensions of the lower-central quotients do not impose a uniform
bound on the degree at which new commutator information may first
appear. More strongly, even the complete underlying lower-central
graded vector space does not impose such a bound.
\end{corollary}

\begin{proof}
For every \(m\geq3\),
Theorem~\ref{thm:graded-vector-space-arbitrary-depth} gives nilpotent
Lie algebras \(L_m\) and \(M_m\) with isomorphic underlying
lower-central graded vector spaces and identical graded commutator
operations through degree \(m-1\), but with nonisomorphic graded Lie
algebra structures detected in degree \(m\). Since \(m\) is
arbitrarily large, neither the dimensions of the lower-central
quotients nor their underlying graded vector space can provide a
uniform bound on the degree at which new commutator structure first
appears.
\end{proof}

\begin{remark}
The point is stronger than equality of lower-central ranks. In the
examples above the entire underlying lower-central graded vector space
is the same; what changes is the graded commutator operation. Thus the
unbounded depth is caused by genuinely higher structural information,
rather than by variation in the sizes of the successive layers.
\end{remark}

Canonical reduction therefore separates the formal length of a lift
tower from its genuinely distinguishing depth. The class-three example
shows that successive commutator stages can each be essential, while
the arbitrary-depth construction shows that essential commutator
information may occur arbitrarily late. More strongly, this phenomenon
persists over the single fixed abelianisation shadow $A=F^4$ and even
after fixing the complete underlying lower-central graded vector
space. Thus classification failure is not merely a consequence of beginning
with a coarse shadow. Hidden structure may remain unresolved through
arbitrarily long sequences of increasingly refined commutator data.
The depth of classification failure is therefore itself a structural
feature of the lifting problem. 

\section{Conclusion}
\label{sec:conclusion}

This paper develops a general framework for studying the hidden
structure omitted by shadow maps. Rather than treating the information
lost under a shadow as an isolated phenomenon, the framework organises
it into lift packages, lift towers, and numerical measures of lift
depth, providing a systematic description of how hidden structure may
be progressively recovered. The starting point is the observation that many mathematical settings
admit natural shadow maps which retain visible information while
omitting finer structural features. The hidden structure is localised
within shadow fibres, where lift packages encode information not
detected by the shadow. Organising this information into successive
layers leads naturally to lift towers, while lift height measures the
depth of hidden structure relative to the chosen shadow.

The framework was developed concretely through nilpotent groups and Lie
algebras. In the quadratic setting, the first hidden structure beyond
abelianisation is encoded by the commutator and quadratic maps. For
groups of higher nilpotency class, the successive lower-central
quotients and their commutator operations form a natural commutator
lift tower whose layers recover the associated graded Lie structure.
For a nilpotent group of class \(c\), this tower has commutator lift
height \(c-1\). A central theme of the paper is that lift height is fundamentally
relative to the chosen shadow. Different shadows make different amounts
of structure visible and therefore leave different amounts of hidden
structure to be recovered. Relative commutator lift height makes this
dependence explicit, while shadow refinement provides a systematic way
to compare the resulting towers. Under the corresponding compatibility
conditions, refining the shadow can only decrease the remaining lift
depth.

The internal structure of lift towers leads to a sharper notion of
depth. Every finite lift tower induces a refinement chain of equivalence
relations on the corresponding shadow fibre. Removing repeated stages
produces a canonical reduced refinement chain, and essential lift
height records only those stages that give strict refinements. The
class-three construction shows concretely that consecutive commutator
lift stages can both be essential, so higher stages need not merely
repeat distinctions already visible at lower levels. More strongly, essential commutator information can occur arbitrarily
deep in a lift tower. For every \(m\geq3\), there exist nilpotent Lie
algebras of class \(m\) whose commutator lift packages agree through
degree \(m-1\) but are distinguished by the degree-\(m\) operation.
Consequently, classification failure can persist through arbitrarily
long finite commutator lift towers even over the single fixed
abelianisation shadow \(A=F^4\). In particular, there is no uniform
finite lift-depth bound determined by the dimension of the visible
shadow.

The obstruction is stronger than a failure of lower-central dimensions
to control classification depth. The arbitrary-depth examples may be
chosen so that their associated lower-central graded objects are
isomorphic as graded vector spaces and their graded commutator
operations agree through arbitrarily high degree, while their graded
Lie algebra structures are distinguished by the next operation. Thus
even the complete underlying lower-central graded vector space does not
bound the degree at which genuinely new commutator information may
appear. Classification depth is therefore controlled not merely by the
sizes of the successive hidden layers, but by the operations linking
those layers. These results sharpen the conceptual role of lift theory. Lift packages
record hidden structure beyond a chosen shadow, lift towers organise
its successive recovery, and reduced lift towers distinguish formal
tower length from genuinely distinguishing depth. The resulting
essential lift height measures how long classification ambiguity can
survive as progressively finer structural information is introduced.
In this sense, the depth of classification failure is itself a
structural feature of the lifting problem.

Although the principal examples considered here arise from nilpotent
groups and Lie algebras, the abstract framework depends only on the
presence of a shadow map together with suitable lift data. Natural
directions for further work include identifying settings in which lift
packages or reduced lift towers arise canonically, developing stronger
functorial constructions, relating lift layers to obstruction and
extension theories, and determining whether analogous unbounded-depth
phenomena occur in other algebraic and categorical classification
problems.

\bibliographystyle{amsplain}
\bibliography{refs}

\end{document}